\documentclass{m2an}

\usepackage[T1]{fontenc}
\usepackage[utf8]{inputenc}
\usepackage{lmodern}
\usepackage[british]{babel}
\usepackage{changes}

\usepackage{amssymb, amsmath, amsthm, amsopn}
\usepackage{thmtools}
\usepackage{mathtools}
\usepackage{commath}
\usepackage{nicefrac}
\usepackage{array}
\usepackage{multirow}

\usepackage{bm}
\usepackage{stmaryrd}
\usepackage{etoolbox}

\numberwithin{equation}{section}

\usepackage{booktabs}
\usepackage{hyperref}

\usepackage{graphicx}
\usepackage{tikz}
\usepackage{pgfplots}
\usepgfplotslibrary{colorbrewer,groupplots}

\usepackage{tikz-3dplot}
\usetikzlibrary{positioning, shapes, arrows, calc, arrows.meta, patterns, shadows, trees, intersections, spath3, calligraphy, spy}

\pgfplotsset{
  compat=1.15,
  cycle list/Set1-5,
  title style = {font = \small},
  legend style = {font = \small},
  label style = {font = \footnotesize},
  tick label style = {font = \footnotesize},
  PPk/.style={mark=square*},
  PPkw/.style={mark=halfcircle*},
  PPkag/.style={mark=10-pointed star},
  TTk/.style={mark=pentagon*},
  TTkag/.style={mark=star},
  TTkemb/.style={mark=triangle*},
  TTkw/.style={mark=halfdiamond*},
}

\usepackage{marginnote}

\newcommand*{\vertbar}{\rule[-1ex]{0.5pt}{8.5ex}}
\newcommand*{\horzbar}{\rule[.5ex]{5.3ex}{0.5pt}}

\newcommand{\Th}{\mathcal{T}_h}

\newcommand{\Fh}{\mathcal{F}_h}
\newcommand{\cC}{\mathcal{C}}

\newcommand{\ThG}{\Th^{\Gamma}}

\newcommand{\Fhgp}{\Fh^\text{gp}}
\newcommand{\Fhgpmin}{\Fh^{\text{gp},\min}}
\newcommand{\Fhgpstar}{\Fh^{\text{gp}\star}}

\newcommand{\Thag}{\Th^{\text{ag}}}

\newcommand{\Ch}{\cC_h}
\newcommand{\ThC}{\Th^\mathcal{C}}
\newcommand{\Thw}{\Th^\omega}
\newcommand{\Tw}{\mathcal{T}_{\omega_F}}
\newcommand{\FhTw}{\Fh(\Thw)}

\newcommand{\Ex}{\mathcal{E}}
\newcommand{\inter}{\pi_h}

\newcommand{\PP}{\mathbb{P}}

\newcommand{\TT}{\mathbb{T}}

\newcommand{\RR}{\mathbb{R}}

\newcommand{\Giso}{\mathcal{G}}
\newcommand{\GisoT}{\mathcal{G}_{\TT}}
\newcommand{\GagT}{\mathcal{G}_{\text{ag}}}

\DeclareMathOperator{\range}{range}
\DeclareMathOperator{\dist}{dist}
\DeclareMathOperator{\id}{id}
\DeclareMathOperator{\Int}{Int}
\DeclareMathOperator{\meas}{meas}

\newcommand{\Bb}{\mathbf{B}}
\newcommand{\Wb}{\mathbf{W}}
\newcommand{\wb}{\mathbf{w}}
\newcommand{\lb}{\mathbf{l}}
\newcommand{\Tb}{\mathbf{T}}
\newcommand{\eb}{\mathbf{e}}
\newcommand{\xbu}{\mathbf{x}}
\newcommand{\ub}{\mathbf{u}}
\newcommand{\ubT}{\ub_{\TT}}
\newcommand{\ubTag}{\ub_{\TT,\text{ag}}}

\newcommand{\vb}{\mathbf{v}}

\newcommand{\orange}[1]{{\color{orange!40!red}#1}}

\newcommand{\BbT}{\Bb_{\TT}}
\newcommand{\lbT}{\lb_{\TT}}

\newcommand{\xb}{\bm{x}}
\newcommand{\xbh}{\hat{\xb}}
\newcommand{\nb}{\bm{n}}
\newcommand{\nbF}{\nb_F}
\newcommand{\nbFh}{\hat{\nb}_F}

\newcommand{\usim}{\tilde{u}}
\newcommand{\uT}{u_{\TT}}

\newcommand{\Ah}{\mathcal{A}_h}
\newcommand{\Bh}{\mathcal{B}_h}

\newcommand{\Sh}{\mathcal{S}_h}
\newcommand{\Lh}{\mathcal{L}_h}

\newcommand{\A}{\mathcal{A}}
\renewcommand{\L}{\mathcal{L}}

\newcommand{\uhom}{u_{\text{hom}}}
\newcommand{\Ahg}{{A}_h}
\newcommand{\Bhg}{{B}_h}
\newcommand{\Shg}{{S}_h}
\newcommand{\Lhg}{{L}_h}
\newcommand{\Whg}{{W}_h}

\newcommand{\nF}{\nb_F}
\newcommand{\dnF}{\nF\cdot\nabla}
\newcommand{\dn}{\nb\cdot\nabla}
\newcommand{\jump}[1]{\llbracket#1\rrbracket} 
\newcommand{\mean}[1]{\{\!\!\{#1\}\!\!\}}
\newcommand{\sigb}{\bm{\sigma}}

\newcommand{\OT}{\Omega_{\mathcal{T}}}

\newcommand{\nrm}[1]{\Vert #1 \Vert}
\newcommand{\snrm}[1]{\vert #1 \vert}
\newcommand{\tripnrm}[1]{\mathopen{|\mkern-1.5mu|\mkern-1.5mu|}#1\mathclose{|\mkern-1.5mu|\mkern-1.5mu|}}

\newcommand{\tnrmA}[1]{\tripnrm{#1}_{\Ah}}
\newcommand{\snrmS}[1]{\vert #1 \vert_{\Sh}}
\newcommand{\tnrmB}[1]{\tripnrm{#1}_{\Bh}}

\newcommand{\tnrmAg}[1]{\tripnrm{#1}_{\Ahg}}
\newcommand{\snrmSg}[1]{\vert #1 \vert_{\Shg}}
\newcommand{\tnrmBg}[1]{\tripnrm{#1}_{\Bhg}}

\DeclareMathOperator{\diam}{diam}
\DeclareMathOperator{\atan}{atan}

\newcommand{\cPhi}{\circ\Phi_h}
\newcommand{\cPhii}{\cPhi^{-1}}
\newcommand{\DPiT}{D\Phi_h^{-T}}
\newcommand{\dDP}{\det(D\Phi_h)}
\newcommand{\ut}{\widetilde{u}}
\newcommand{\wht}{\widetilde{w}_h}

\newcommand{\eqrefh}[1]{\textup{(\hyperref[#1]{\ref*{#1}\textsuperscript{hom}})}}
\newcommand{\ghom}{g^\textup{hom}}

\theoremstyle{definition}
\newtheorem{ssmptn}[thrm]{Assumption}

\begin{document}

\title{Unfitted Trefftz discontinuous Galerkin methods for elliptic boundary value problems}

\thanks{HvW acknowledges support through the Austrian Science Fund (FWF) project F65.}

\author{Fabian Heimann}\address{Institut für Numerische und Angewandte Mathematik, University of
  Göttingen, Lotzestraße 16-18,
  37083 Göttingen, Germany; \email{ \{f.heimann, lehrenfeld, p.stocker\}@math.uni-goettingen.de}}
\author{Christoph Lehrenfeld}\sameaddress{1}
\author{Paul Stocker}\sameaddress{1}
\author{Henry von Wahl}\address{Fakultät für Mathematik, Universität Wien, Oskar-Morgenstern-Platz, 1090 Wien, Austria; \email{henry.wahl@univie.ac.at}}

\date{\today}

\begin{abstract}
  We propose a new geometrically unfitted finite element method based on discontinuous Trefftz ansatz spaces. 
Trefftz methods allow for a reduction in the number of degrees of freedom in discontinuous Galerkin methods, thereby, the costs for solving arising linear systems significantly. 
  This work shows that they are also an excellent way to reduce the number of degrees of freedom in an unfitted setting.
  We present a unified analysis of a class of geometrically unfitted discontinuous Galerkin methods with different stabilisation mechanisms to deal with small cuts between the geometry and the mesh. We cover stability and derive a-priori error bounds, including errors arising from geometry approximation for the class of discretisations for a model Poisson problem in a unified manner. The analysis covers Trefftz and full polynomial ansatz spaces, alike. Numerical examples validate the theoretical findings and demonstrate the potential of the approach.
\end{abstract}

\maketitle

\section{Introduction}
In the last two decades, unfitted finite element methods became popular as an 
alternative to more traditional body-fitted methods to solve
partial differential equations on complex geometries numerically. The idea is to separate
the computational mesh from the geometry description to remove the burden of
mesh generation, mesh adaptation and remeshing when dealing with complex (and
possibly time-dependent) geometries. 
This is accomplished by embedding the domain of interest into an unfitted background mesh.
In the context of finite element methods,
unfitted discretisations go under different names such as 
\emph{CutFEM}~\cite{B15},
\emph{extended FEM} (X-FEM)\cite{BMUP01,MDB99}, \emph{Finite Cell} \cite{PDR07} and several more.
In many cases, discontinuous Galerkin (DG) methods are attractive on unfitted meshes as they are
for fitted meshes, e.g., for convection-dominated convection-diffusion problems, because of their flexibility in regards to the polynomial basis function as exemplified by Trefftz and $hp$-methods or computational aspects.
In the spirit of CutFEM and X-FEM, there are several unfitted discretisations based on a
discontinuous Galerkin formulation, cf.~\cite{BE09,EH12,Mas12,MKO13,BHLM17,K17,GM19}.

\medskip

A challenge for unfitted methods arises from the fact that the boundary can cut through mesh elements arbitrarily. 
Trimming parts of the mesh outside the domain of interest may lead to ill-shaped elements.
Two mechanisms have proven fruitful in making such methods robust with respect to ill-shaped elements: \emph{ghost penalties} and \emph{element aggregation}. 

Element aggregation (\texttt{AG}) 
joins mesh elements to ensure that the support of basis functions does not degenerate for bad-cut configurations in the mesh. 
This has been done among others in \cite{EH12,K17} for DG methods and in \cite{BCDE21} for unfitted Hybrid-High-Order methods.
In \cite{BVM18} and the proceeding works of Verdugo and Badia, the idea has also been generalised to continuous finite elements. 
The method is also referred to as cell merging or cell agglomeration.

Ghost penalty (\texttt{GP}) 
introduces an additional stabilisation term, the
ghost penalty stabilisation \cite{Bur10,B15}, that introduces a
\emph{volumetric} coupling, i.e., a coupling that involves all unknowns of two
adjacent elements in the vicinity of shape-irregular cuts.
No mesh elements or basis functions are changed for this approach.
In order to reduce the number of coupled elements, the stabilisation term can be applied within patches only. These patches are built by the same machinery used by element aggregation, see \cite{BNV22}.
We refer to this method as patch-wise ghost penalty ($\omega$\texttt{GP}); it also appears under the names of weak ghost penalty or weak aggregated elements. We also mention the work \cite{ELL18} for an approach where ghost penalties are avoided. The ill-cut problem is resolved here by removing basis functions whose supports have small intersections with the computational domain.

\medskip

Discontinuous Galerkin methods come with increased 
degrees of freedom (dofs) when compared to their continuous counterpart.
When solving linear systems for corresponding discretisations, the duplication
of degrees of freedom affects the efficiency of numerical methods even more.
There are two well-known remedies in the literature for the body-fitted case.

The first is to use Hybrid Discontinuous Galerkin (HDG) methods \cite{CGL09}, where
additional degrees of freedom are introduced on element interfaces. These
additional unknowns allow (in many cases) the elimination of interior (volumetric)
DG unknowns by a Schur complement strategy (known as static condensation). The
remaining degrees of freedom in the global linear system are then
significantly reduced, especially in the case of
higher-order discretisations. A difficulty with the extension of HDG methods
from the fitted to the unfitted case is the robustness with respect to
shape-irregular cuts. Applying a ghost penalty stabilisation is not possible as
the couplings between direct element neighbours introduced by the stabilisation
terms contradict the decoupling exploited in the hybridisation. Cell merging
strategies are possible but require the handling of polygonal meshes with
corresponding facet functions as it is done in the unfitted Hybrid High Order (HHO)
method \cite{BE17,BE18}. In \cite{GKF17}, neither ghost penalties nor cell
merging needs to be used by changing the background mesh to avoid ill-shaped cuts
(in two space dimensions).

The second remedy to overcome the computational costs of DG methods is the
class of Trefftz DG methods. In Trefftz methods, originating from 
\cite{trefftz1926}, the ansatz space is constructed to lie in the kernel of the differential operator of the PDE at hand.
Compared with a corresponding DG method, the same accuracy can be
achieved for these methods at significantly reduced costs. 
Trefftz-DG methods for the Laplace problem are analysed in \cite{HMPS14, LiShu2012, LiShu2006}.
Indeed, the complexity reduction is comparable to that of the HDG method, cf.~\cite{LS22}.
In contrast to the HDG mechanism, the mechanism to reduce the computational
complexity does not interfere with either the ghost penalty stabilisation or
the cell merging strategy. 

\subsection{Main Contributions and Outline of the Paper} \label{sec:maincontr}

In this work, we will consider the
combination of unfitted discontinuous Galerkin formulations with Trefftz DG
finite element spaces on the example of the Poisson problem.
Our objective is to develop the tools for analyzing and advancing the method beyond this particular model problem.
To the best of our knowledge, this is the first occasion
of a geometrically unfitted Trefftz DG method in the literature.

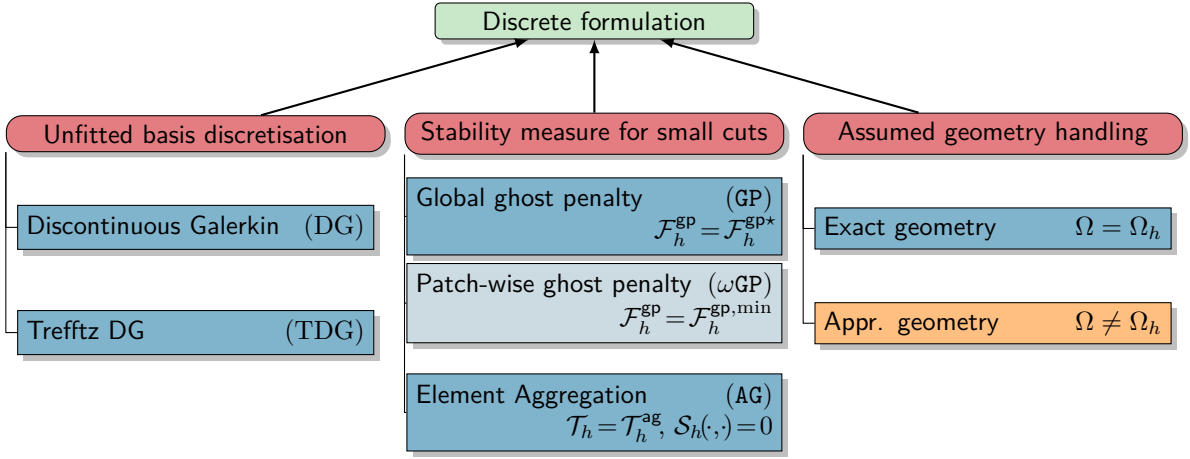
\begin{figure}
  \centering
  \tikzset{
  basic/.style  = {draw, text width=6cm, drop shadow, font=\sffamily, rectangle},
  root/.style   = {basic, text width=4cm, rounded corners=2pt, thin, align=center,
    fill=green!40!gray!30},
  level 2/.style = {basic, text width=4.8cm, rounded corners=6pt, thin,align=center, fill=orange!30!purple!60},
  level 3/.style = {basic, text width=4.5cm, thin, align=left, fill=blue!60!green!50}
}

\begin{tikzpicture}[
    level 1/.style={sibling distance=5.3cm},edge from parent/.style={<-,draw,thick},>=latex]
  \tikzstyle{every node}=[font=\small]
\node[root] {Discrete formulation}
child {node[level 2] (c1) {Unfitted basis discretisation}}
  child {node[text width=5.5cm, level 2] (c2) {Stability measure for small cuts}}
  child {node[level 2] (c3) {Assumed geometry handling}};

  \begin{scope}[every node/.style={level 3}]
    \node [below =0.7cm of  c1] (c11) {Discontinuous Galerkin  \hfill \eqref{eqn.poisson-poblem.dg}};
      \node [below =0.8cm of c11] (c12) {Trefftz DG \hfill \eqref{eqn.poisson-poblem.trefftz}};

    \node [text width=4.75cm,below =0.3cm of c2] (c21) {Global ghost penalty \hfill \eqref{case:gp} \\ \hfill $\Fhgp\!=\!\Fhgpstar\!$};
    \node [text width=4.75cm,below =0.1cm of c21, fill=blue!60!green!70!black!20] (c22) {Patch-wise ghost penalty \hfill \eqref{case:pgp} \\ \hfill $\Fhgp\!=\!\Fhgpmin$};
    \node [text width=4.75cm,below =0.4cm of c22] (c23) {Element Aggregation \hfill \eqref{case:ag}\\ \hfill $\Th\!=\!\Thag\!,\,\Sh\!(\!\cdot,\!\cdot\!)\!=\!0$};

    \node [below =0.7cm of c3] (c31) {Exact geometry \hfill $\Omega = \Omega_h$};
    \node [below =0.7cm of c31, fill=orange!50] (c32) {Appr. geometry \hfill $\Omega \neq \Omega_h$};
  \end{scope}

  \foreach \value in {1,2}
  \draw[-] (c1.south west) |- (c1\value.west);

  \foreach \value in {1,2,3}
  \draw[-] (c2.south west) |- (c2\value.west);

  \foreach \value in {1,2}
  \draw[-] (c3.south west) |- (c3\value.west);  
\end{tikzpicture}     \caption{Overview on different discretisation settings treated in this manuscript. Twelve different settings arise from the possible choices in the basis discretisation (DG vs. Trefftz DG), stability measure for small cuts (ghost penalty in two versions vs. aggregation) and the assumption on the geometry handling (exact vs. approximated).}
\label{fig.method-zoo}
\end{figure}

Our analysis covers the unfitted DG method with two choices for the basis functions: either the full polynomial space or the Trefftz space. 
These choices of basis functions are then combined with either the cell-merging, the ghost penalty, or the patch-wise ghost penalty stabilisation -- cf.\@{} also Figure~\ref{fig.flow-diagramm} below for a comparison of stabilisation strategies --  to arrive at robust unfitted discretisations, which are then analysed. 
The analysis covers a higher-order a-priori error analysis, which we first present for exact geometries and then extend to the case with geometry approximation errors.
A summary of these methods presented here is given in Figure~\ref{fig.method-zoo}.

The unfitted DG method with ghost penalty has already been proposed and 
analysed in~\cite{GM19}, except for the analysis of the geometrical errors. 
The element aggregation and patch-wise ghost penalty follow from the 
works~\cite{BNV22} and~\cite{BVM18} on the (more general) aggregated FEM.
Our unified analysis is possible, as both the Trefftz and aggregated finite 
element spaces are subsets of the standard DG space. 
The patch-wise ghost penalty consists of a minor modification of the 
ghost-penalty term.

Novel contributions of the work are the
\begin{itemize}
  \item description of unfitted DG and unfitted Trefftz DG methods with three different stabilisation mechanisms, both with and without geometry approximation errors,
  \item unified analysis leading to a priori error estimates for these unfitted (Trefftz) DG methods, including geometry errors, and
\item numerical experiments for these methods and the discussion of implementational aspects.
\end{itemize}

The remainder of this paper is structured as follows: In Section~\ref{sec.intro.subsec.problem}, we present the model problem under consideration in this paper. In Section~\ref{sec.DG-TDG-exact-geometry}, we recap the unfitted DG method under the assumption of exact geometry handling, as covered in previous literature, and introduce the Trefftz DG method under the same assumption on the geometry. 
In Section~\ref{sec.stabilization}, we consider the different approaches to deal with the problem of small cuts, namely element aggregation in Section~\ref{sec.stabilization:subset.aggregation} and ghost-penalty stabilisation in Section~\ref{sec.stabilization:subset.ghost-penalty}.
Section~\ref{sec.error-exact-geometry} presents the unified error analysis for the considered methods. The analysis extends the work on unfitted DG methods in \cite{GM19}. A crucial role for this extension is the special choice for the interpolant of smooth functions.
Section~\ref{sec.anaysis-with-geom-error} then covers the error analysis of the unfitted DG and TDG methods, including geometry approximation errors inherent in unfitted finite element methods.
We then discuss several variants and implementational aspects of the covered methods in Section~\ref{sec.further-methods}, including embedding the Trefftz and aggregated spaces into standard discontinuous finite element spaces.
We present numerical examples of the methods in Section~\ref{sec.numerical-examples}. Here we include examples with and without geometry approximation errors. Finally, in Appendix~\ref{sec.appendix}, we give some proofs for completeness that consist only of minor adaptations of proofs available in the literature.

\subsection{Model Problem}
\label{sec.intro.subsec.problem}
As a model boundary value problem, we consider the Poisson problem on an open bounded domain $\Omega \subseteq \mathbb{R}^d, d = 2,3$ with boundary $\Gamma \coloneqq \partial \Omega \in \xCtwo$: Let $f \in \xLtwo(\Omega)$ and $g \in \xHn{\frac{1}{2}}(\Gamma)$ be given. Then the problem reads: Find $u\colon \Omega \to \mathbb{R}$ such that
\begin{subequations}\label{eqn.strong-poisson}
  \begin{align}
 -\Delta u &= f \quad \textnormal{in } \Omega \label{eqn.strong-poisson.bulk}\\
 u &= g \quad \textnormal{on } \Gamma.\label{eqn.strong-poisson.boundary}
\end{align}
\end{subequations}
Defining $V_g \coloneqq \{ v \in \xHn{1}(\Omega) \, | \, v|_\Gamma = g \}$, we can state the following weak form to the strong form given above: Find $u \in V_g$ such that
\begin{equation*}
 (\nabla u,  \nabla v)_\Omega = (f,v)_\Omega \quad \forall v \in V_0.
\end{equation*}
Here, we used the inner product $(u,v)_\Omega \coloneqq \int_{\Omega} uv \mathrm{d}x$. We note that for the analysis below, we will homogenise the problem with respect to the volume forcing, meaning that we compute a particular solution of \eqref{eqn.strong-poisson.bulk} and then discretise an homogeneous problem, see \hyperref[eqn.strong-poisson.hom]{(\ref*{eqn.strong-poisson}\textsuperscript{hom})} below. This is standard for Trefftz methods. We will also discuss the numerical realisation of this homogenisation in Section~\ref{sec.further-methods} and present a numerical example with $f\neq 0$.

\section{Unfitted DG and Trefftz DG Methods With Exact Geometry}
\label{sec.DG-TDG-exact-geometry}
In this section, we first recap the unfitted DG method assuming exact geometries, i.e.\@{} neglecting 
errors arising from inaccurate integration over cut elements.
Based on this, we will then introduce the unfitted Trefftz DG method.

\subsection{Preliminaries: Geometry, Mesh and Cut Elements}\label{sec:cutel} 
We start by introducing some notation and assumptions: Let $\tilde \Omega$ be a
background domain, sufficiently large such that $\overline{\Omega} \subseteq
\tilde \Omega$ and let $\widetilde{\mathcal{T}}_h = \{T\}$ be a division of $\tilde\Omega$ into non-overlapping shape-regular elements. 
The local mesh size of a mesh element $T\in\widetilde{\mathcal{T}}_h$ is
defined as $h_T = \diam(T) \coloneqq \sup_{\xb_1, \xb_2 \in T} \nrm{\xb_1 -
\xb_2}_2$. We allow for quite general polyhedral meshes including possibly curved elements,
but require that the mesh conforms to the following assumption or its subsequent
relaxed version.

\begin{ssmptn}\label{assumption.mesh1}
  We assume that for every $T' \in \widetilde{\mathcal{T}}_h$ there holds
\begin{itemize}
  \item[(a)] There are two balls $b_{T'} \subset {T'} \subset B_{T'}$, such that ${T'}$ is star shaped with respect to the ball $b_{T'}$ and $\diam(B_{T'})/\diam(b_{T'}) \lesssim 1$.
  \item[(b)] The element boundary can be divided into mutually exclusive subsets $\{ F_i \}_{i=0}^{n_{T'}}$ with $\diam(T') \leq c \diam(F_i)$, $i=0,...,n_{T'}$, where $n_{T'}$ and $c$ are uniformly bounded, satisfying
  \begin{itemize}
    \item[(i)] There exists a sub-element ${T'}_{F_i}$ of $T'$ with $d$ planar facets meeting at a
    vertex $\xb^0_i\in T'$, such that  ${T'}_{F_i}$ is star-shaped with respect to $\xb^0_i$ and $h_{T_{F_i}'} \simeq h_{T'}$.
    \item[(ii)] There exists a uniform constant $c_{\ref{assumption.mesh1}}$, such that 
    \begin{equation*}
      (\xb - \xb^0_i) \cdot \nb_{F_i}(x) \geq c_{\ref{assumption.mesh1}} h_{T'} \quad \forall \xb \in F_i.
    \end{equation*}
  \end{itemize}
  \item[(c)] The element boundary $\partial {T'}$ is the union of a finite (yet, arbitrarily large) number of closed $C^1$ surfaces.
\end{itemize}
\end{ssmptn}
Here and in what follows, we use the notation $a \lesssim b$ if there exists a constant $c>0$, independent of the mesh size and mesh-interface cut position, such that $a\leq cb$. Similarly, we use $\gtrsim$ if $a\geq cb$, and $a\simeq b$ if both $a\lesssim b$ and $b \lesssim a$ holds. These assumptions are essentially based on \cite[Assumptions 4.1 and 4.3]{CDG21} to guarantee that the trace and inverse estimates cited from this work are valid here. For further details on these mesh assumptions, we refer to \cite[Section~4, Figure~2 and Figure~3]{CDG21}. For a simpler but more restrictive mesh assumption for polytopic meshes under which appropriate trace estimates are available, we also refer to \cite{CDGH17}.

Let us note that (possibly curved) simplicial, hexahedral or quadrilateral meshes that are 
shape-regular in the usual sense fulfil Assumption~\Rref{assumption.mesh1}.
In the following, we want to enlarge the class of admissible meshes to those
arising from merging a (uniformly bounded) number of (neighbouring) elements from
meshes fulfilling Assumption~\Rref{assumption.mesh1}. A corresponding relaxed
version of Assumption~\Rref{assumption.mesh1} is the following.
\begin{ssmptn}\label{assumption.mesh2}
  For every $T\in\widetilde{\mathcal{T}}_h$, we assume that $T$ is a Lipschitz
  domain and that it can be decomposed in $m_T$ non-overlapping elements
  $\{T'\}$, with $m_T$ uniformly bounded and for each element $T'$ the
  assumptions in Assumption~\Rref{assumption.mesh1} holds.
\end{ssmptn}

\begin{lmm} \label{lmm.ass1ass2}
  Assumption~\Rref{assumption.mesh1} directly implies
  Assumption~\Rref{assumption.mesh2} with $m_T=1$ and $\{T'\} = \{T\}$ for every
  $T \in \widetilde{\mathcal{T}}_h$. Further, it immediately follows that for every
  sub-element $T'$ of an element $T\in\widetilde{\mathcal{T}}_h$ with
  $\widetilde{\mathcal{T}}_h$ fulfilling Assumption~\Rref{assumption.mesh2}  there are two balls $b_{T'}
  \subset T' \subset T \subset B_{T}$ such that $T'$ is star shaped w.r.t. the
  ball $b_{T'}$ and $\diam(B_{T})/\diam(b_{T'}) \lesssim m_T$.
\end{lmm}

\begin{rmrk}\label{rmrk.assumptions}
  The first part of Assumption~\Rref{assumption.mesh1}
  (even in its relaxed version of Assumption~\Rref{assumption.mesh2}) guarantees
  a shape regularity property sufficient for an inverse estimate and
  the interpolation via Taylor polynomials below, as needed for the analysis of
  the Trefftz method. 
  The second assumption is essentially
  a bound on the curvature of the elements, although this
  assumption is a weak restriction compared with
  Assumption~\Rref{assumption.mesh1}(a).
  By construction, starting from a mesh fulfilling Assumption~\Rref{assumption.mesh1} and applying
  a cell merging strategy (where always only a uniformly bounded number of neighbouring
  elements are merged) results in a mesh fulfilling
  Assumption~\Rref{assumption.mesh2}. Obviously, further merging of a resulting
  mesh only fulfilling Assumption~\Rref{assumption.mesh2} will still yield a mesh
  fulfilling Assumption~\Rref{assumption.mesh2}. However, these merging
  steps will decrease the shape regularity bound by a multiplicative constant ($\sim m_T^{-1}$).
\end{rmrk}

Of specific interest will be the following parts of the background mesh, which we call the active mesh and the cut mesh, i.e.\@{} the parts of the background mesh that contribute to the covering of $\Omega$ and the part that is intersected by the domain boundary $\Gamma$, respectively:
\begin{equation*}
  \Th = \{ T \in \widetilde{\mathcal{T}}_h \::\: \meas_{d}(T \cap \Omega)>0 \},
  \quad
  \ThG  = \{ T \in \widetilde{\mathcal{T}}_h \::\: \meas_{d-1}(T \cap \Gamma)>0 \}.
\end{equation*}
We collect the domain of the active mesh as $\OT \coloneqq \Int\bigcup_{T\in\Th}\overline{T}$.

We further introduce sets of facets needed for the unfitted DG method. For the communication between all direct neighbour elements in a set of elements $S = \{T\}$ we collect
\begin{equation}\label{eqn.sets-of-facets}
 \Fh(S) = \{ F = \partial T_1 \cap \partial T_2 \::\: T_1, T_2 \in S, T_1 \neq T_2 \text{ and } \meas_{d-1}(F)>0 \},
\end{equation}
and denote $\Fh = \Fh(\Th)$.

With abuse of notation we denote by $h$ the global mesh size $h = \max_{T\in \Th} h_T$ when used as a scalar, as well as the piecewise constant field on $\Omega$, $h: \Omega \to \mathbb{R}$ with $h|_T = h_T$, $T \in \Th$ or as the piecewise constant field on the skeleton, $h: \Fh \to \mathbb{R}$ with $h|_F = h_F$, $F \in \Fh$ with $h_F = \diam(F) = \sup_{\xb_1, \xb_2 \in F} \nrm{\xb_1 -
\xb_2}_2$. Note that due to shape regularity $h_F \leq h_T \lesssim h_F$ for $F \subset \partial T$.

\subsection{Unfitted DG Methods With Exact Geometry} 
Starting point and first method is the unfitted DG discretisation as in \cite{GM19}.
The discrete function spaces are given as the discontinuous polynomials of order $k$ on $\Th$:
\begin{equation}\label{eqn.dg-space}
 \PP^k(\Th) = \bigoplus_{T \in \Th} \{ p \in \PP^k(T) \},
\end{equation}
where $\PP^k(T)$ is the space of polynomials up to degree $k$ on $T$.

As usual with interior penalty DG methods, we penalise jumps across facets $F \in \Fh$. For this, we need the following average and jump operations:
\begin{align*}
 \mean{\sigma }|_F &= \frac{1}{2} (\sigma_F^+ + \sigma_F^-), \quad \mean{\nF \cdot \sigb }|_F = \frac{1}{2} \nF \cdot  (\sigb_F^+ + \sigb_F^-), \\
 \jump{ w }|_F &= w_F^+ - w_F^-, \quad \textnormal{where} \quad v_F^{\pm}(x) \coloneqq \lim_{t \to 0} v_F(x \pm t \nF),
\end{align*}
and $\nF$ denotes some fixed facet normal to $F$.

Next, in preparation of the discrete variational formulation, we introduce the bilinear form $\Ah$ as
\begin{equation} \label{eq:Ah} \tag{$\Ah$}
\begin{aligned} 
    \Ah(u,v) \coloneqq &\ (\nabla u, \nabla v)_\Omega - (\dn u, v)_\Gamma - (u, \dn v)_\Gamma + \beta (h^{-1} u,v)_\Gamma \\
                       &- (\mean{\dnF u }, \jump{ v })_{\Fh \cap \Omega} - (\jump{u}, \mean{\dnF v } )_{\Fh \cap \Omega} + \beta (h^{-1} \jump{u}, \jump{v})_{\Fh \cap \Omega},
\end{aligned}
\end{equation}
the corresponding linear form $\Lh$ for the right-hand side as
\begin{equation*}
 \Lh(f,g;v) \coloneqq (f,v)_\Omega - (\dn v,g)_\Gamma + \beta (h^{-1} g,v)_\Gamma,
\end{equation*}
and a ghost penalty stabilisation form $\Sh(\cdot,\cdot)$, which we will specify in Section~\ref{sec.stabilization:subset.ghost-penalty}. 

In terms of these definitions, we can then introduce the discrete problem as follows: Find $u_h \in \PP^k(\Th)$ such that
\begin{gather}\label{eqn.poisson-poblem.dg} \tag{DG}
 \Bh(u_h, v_h) \coloneqq \Ah(u_h,v_h) + \Sh(u_h,v_h) = \Lh(f,g;v_h) \quad \forall v_h \in \PP^k(\Th).
\end{gather}

\subsection{Unfitted Trefftz DG Methods With Exact Geometry}
\label{sec.DG-TDG-exact-geometry.subsec.tdg}

One major drawback of discontinuous Galerkin methods is the high computational cost related to the additional degrees of freedom due to the discontinuities. A popular remedy for this is to consider hybridised discontinuous Galerkin (HDG) methods. In HDG methods, additional unknowns are introduced at the facets of elements in order to remove direct couplings between neighbouring element unknowns.
Thereby the size of the linear systems to be solved can be reduced by static condensation significantly, especially for higher order. This approach does not fit well with unfitted finite element methods in general, as stabilisation techniques, such as the ghost penalty method, rely on direct couplings between neighbouring elements and can not be hybridised efficiently. Only if ghost penalty methods can be circumvented, e.g. by the cell aggregation techniques discussed above, a hybridised approach can be applied as in the unfitted HHO methods, see, e.g., \cite{BCDE21}.

An alternative approach to reduce the size of the system to be solved is to consider Trefftz DG methods. Here the essential idea is to use a DG space consisting only of polynomials\footnote{although this can also be generalised to non-polynomial spaces} which fulfil the homogeneous PDE problem on every element.
This leads to a complexity reduction of the number of degrees of freedom as well as the global couplings, which is similar to the HDG method, see the discussion in \cite{LS22}. The construction of the Trefftz DG method does not interfere with the usual coupling pattern of DG methods and can hence be combined with the ghost penalty method straightforwardly.

Let us now introduce a Trefftz version of the previous unfitted DG method simply by replacing the discrete function space to discontinuous polynomials of order $k$, which satisfy the Trefftz condition of the Laplace problem:
\begin{equation}\label{eqn.trefftz-space}
  \TT^k(\Th) = \bigoplus_{T \in \Th} \{ p \in \PP^k(T) \, | \, \Delta p = 0 \}.
\end{equation}
In order to make sense of this subspace of the previous DG space with respect to
the \emph{inhomogeneous} Poisson equation, we assume for now that an
element-wise smooth particular solution $u_p$, with $\Delta u_p = f$ on each element
(possibly discontinuous across element interfaces), is given.

Several approaches exist to homogenise the Laplace problem to apply Trefftz methods; see \cite{zbMATH01240671, zbMATH02139293, LS22, yang2020trefftz}.
The approaches \cite{zbMATH01240671, zbMATH02139293, yang2020trefftz} focus on collocation-based methods.
The embedded Trefftz method, presented in \cite{LS22}, is DG based and gives an easy way to homogenise the system, which we review in Section~\ref{sec.further-methods}.

In terms of these definitions, we can then introduce the discrete problem as
follows: Find $\uT \in \TT^k(\Th)$ such that $ \forall v_h \in \TT^k(\Th)$ there holds
\begin{equation}\label{eqn.poisson-poblem.trefftz} \tag{TDG}
  \Bh(\uT, v_h) = \Ah(\uT,v_h) + \Sh(\uT,v_h) = \Lh(f, g; v_h) - \Bh(u_p, v_h).
\end{equation}
\begin{rmrk}
The bilinear form $\Ah(\cdot,\cdot)$ with Trefftz test and trial spaces was analysed in \cite{HMPS14} for the fitted case. Using integration by parts again on the first term in $\Ah(\cdot,\cdot)$, one obtains an `ultra-weak formulation', where the volume term vanishes due to the properties of Trefftz test functions, then $\Ah(\cdot,\cdot)$ only requires integration over facets and no additional terms are needed. In the discrete setting, this provides an equivalent formulation and can bring considerable savings for the assembly of the linear system. For the Helmholtz problem, such an ultra-weak formulation has been studied in \cite{CD98,HMP16}.
\end{rmrk}

\section{Stabilisation Techniques}
\label{sec.stabilization}

In this section, we consider different approaches to deal with stability in the presence of shape-irregular cut configurations. These approaches are element aggregation and different ghost penalty stabilisations.
In the remainder of this work, we will cover all variants in a mostly unified manner.

\subsection{Element Aggregation}
\label{sec.stabilization:subset.aggregation}

We can avoid the presence of shape irregular cut elements by cell-merging strategies. This has been done, among others, in \cite{EH12,JL12,K17} for DG methods and in \cite{BCDE21} for unfitted Hybrid-High-Order methods. In \cite{BVM18} and the proceeding works of Verdugo and Badia, the idea has also been generalised to continuous finite elements.

We will use the clustering strategies which group certain elements $T
\in \Th$ together in patches $\Thw = \{T\}$, as, e.g., presented in \cite{BVM18}.
These patches are directly used in a cell-merging strategy resulting in an
\emph{aggregated} element $\omega = \Int(\bigcup_{T\in\Thw} \overline{T})$
per patch. However, they will also prove useful for the strategy based on ghost penalties; see also
\cite{BNV22}.

We denote $\Ch = \{ \omega \}$ as the set of aggregated elements obtained from the aggregation of two or more elements, and $\ThC=\{\Thw,\ \omega\in\Ch\}$ as the set of non-trivially aggregated elements. To every aggregated element we define the diameter $h_\omega = \diam(\omega)$.
Finally, $\Thag$ denotes the active mesh after aggregation. Note that an aggregated element,
i.e.\@{} an element in $\Thag$, may originate only from a single element in the
interior of the domain. In the next paragraph, we summarise the crucial
properties of these aggregated elements.

\begin{figure}\centering
  \includegraphics{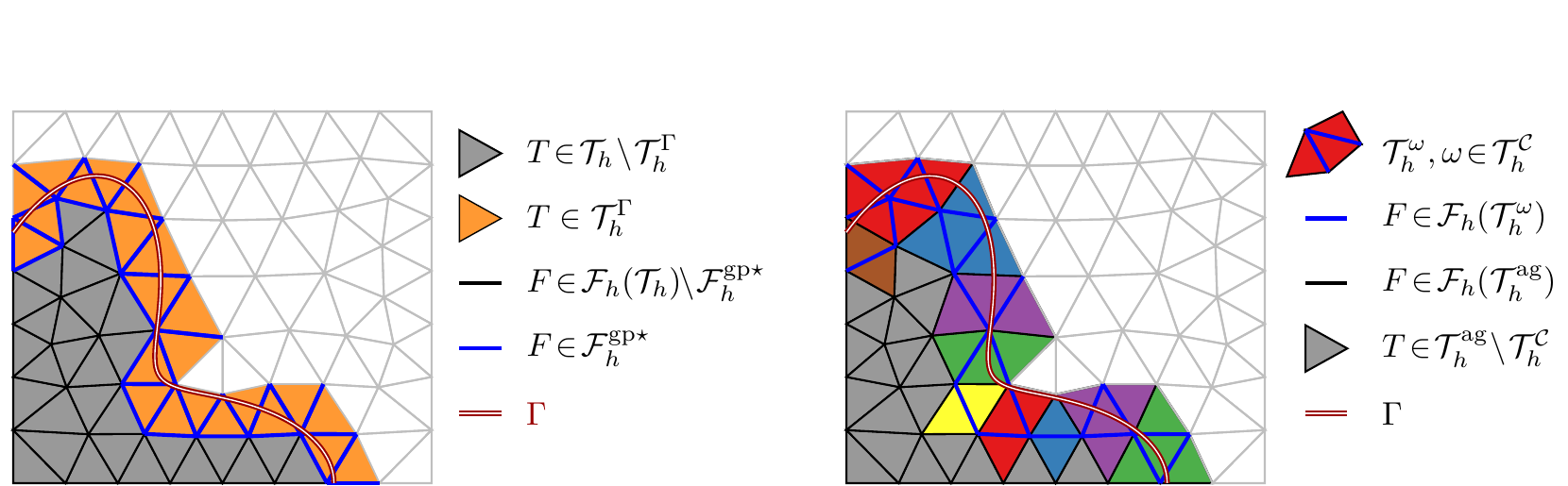}
  \caption{Sketches of the different sets of elements and facets for the unfitted methods. Left: Global ghost-penalty stabilisation elements and facets. Right: Aggregated elements (elements gathered in a patch have the same colour) and corresponding inner-patch facets corresponding to patch-wise ghost-penalty stabilisation.}
  \label{fig:discretegeom}
\end{figure}

The purpose of the non-trivial patches is to group every shape-irregular cut element together with at
least one shape-regular \emph{root} element, cf.\@{} Figure~\ref{fig:discretegeom}
for a sketch. Elements that are not adjacent to cut elements are not affected
and form trivial patches. The number of elements in each patch is uniformly
bounded so that the set of aggregated elements $\Thag$ itself fulfils Assumption~\Rref{assumption.mesh2}.

For the construction of the patches we formulate the following assumption. The
fulfilment of this assumption is the (achieved) goal of the strategies in \cite{BVM18}.
\begin{ssmptn}\label{assumption.path-path}
  We assume that from any cut cell $T_0\in\ThG$ (which is inside a patch $\Thw \in \ThC$) there is a path of elements $\{T_0,T_1,\dots, T_n\} \subset \Thw$ such that
  \begin{enumerate}
  \item For $F=\partial T_i \cap \partial T_{i+1}$, we have that $\meas_{d-1}(F\cap\Omega)> 0$,
  \item $T_n\in \Th\setminus\ThG$, and
  \item $n \leq n_\text{max}$, where $n_\text{max}$ is a fixed integer.
  \end{enumerate}
\end{ssmptn}

\begin{rmrk}[Shape-regular vs. shape-irregular elements]
 We note that the construction of the patches in the previous assumption is
 based on a distinction between cut and uncut elements rather than shape-regular
 and shape-irregular cut elements. This could be generalised further,
 e.g., based on $\vert T\cap\Omega \vert / \vert T \vert,~T\in\Th$.
 As it is standard in the CutFEM literature and avoids the introduction of
 additional notational burden, we stay with the simpler choice for ease of presentation.
\end{rmrk}

\begin{lmm}\label{lemma.patch-size}
  Under Assumption~\Rref{assumption.path-path}, the element aggregates
  $\omega=\Int(\overline{T}_0\cup\dots\cup \overline{T}_n)$ have a maximum diameter of $h_\omega \leq (2
  n_\text{max} + 1) \max_{i \in \{0,..,n\}} \{ h_{T_i} \}$. Furthermore, the intersection of the aggregated element with $\Omega$
  is shape regular in the sense that there is a constant $c > 0 $ uniform in
  $\omega \in \Ch$ so that 
  $ | \omega \cap \Omega | \geq c/n | \omega | $.
\end{lmm}

\begin{proof}
  The first statement is stated and proven in \cite[Lemma 2.2]{BVM18}. The
  second follows directly from quasi uniformity.
\end{proof}
We note that the lemma implies that for all $\omega \in \Ch$ we have $h_\omega \simeq h_T$ for all $T \in \Thw$.

After (proper) cell merging, an aggregated mesh is guaranteed to be \emph{cut-shape-regular} in
the following sense:
\begin{dfntn}\label{def.shaperegcut}
  An active mesh $\Th$ that fulfils Assumption~\Rref{assumption.mesh2} is denoted as \emph{cut-shape-regular} if there is a constant
  $c > 0$ (independent of $h$) so that $\min_{T \in \Th} {|T \cap \Omega|}/{|T|} > c$.
\end{dfntn}

\begin{lmm}\label{lemma.shape-regular-norm-equivalence}
For a cut-shape-regular mesh $\Th$, it holds for $v_h\in\PP^k(\Th)$ that
\begin{equation}\label{eq.norm-equivalence.shapereg}
  \nrm{\nabla v_h}_{\OT} \simeq \nrm{\nabla v_h}_{\Omega}.
\end{equation}
\end{lmm}
\begin{proof}
This is an immediate consequence of Definition~\ref{def.shaperegcut}, and the constants in \eqref{eq.norm-equivalence.shapereg} only depend on the cut-shape-regularity constant.
\end{proof}

\subsection{Ghost Penalty and Patch-wise Ghost Penalty}
\label{sec.stabilization:subset.ghost-penalty}

In this section, we assume that we do not have a \emph{cut-shape-regular} mesh and
introduce a stabilisation (in two variants). In this case, we will still use
the patches from the previous section to define regions where the
stabilisation needs to act but not to change the mesh.

Using \eqref{eqn.sets-of-facets}, we define the set of facets $\FhTw$.
This is the set of all interior facets of a patch $\Thw\in\ThC$ and is only needed for the variant \eqref{case:pgp}.

For the ghost-penalty stabilisation, we require a set of facets that connects (possibly indirectly over several elements and facets) every shape-irregular cut element with a shape-regular root element.
We denote this set as $\Fhgp$. As a minimal choice for stability, we take the set of all \emph{interior} facets of all patches
\begin{equation*}
\Fhgpmin \coloneqq \bigcup_{\omega \in \Ch} \FhTw.
\end{equation*}

In this case, the connection from shape irregular to shape regular element is dealt with within each patch separately. A larger set, more often used in the literature, takes all facets between cut and uncut elements
\begin{equation*}
  \Fhgpstar \coloneqq \{ F = \partial T_1 \cap \partial T_2 \, \::\: \, T_1 \in \Th, T_2 \in \ThG, T_1 \neq T_2, \text{ and } \meas_{d-1}(F)>0 \}.
\end{equation*}
For the ghost penalty method, we set $\Fhgp = \Fhgpstar$; for the patch-wise ghost penalty method, we set $\Fhgp = \Fhgpmin$. 
In the analysis, we will need $\Fhgp \supset \Fhgpmin$ to prove stability and $\Fhgp \subset \Fhgpstar$ to prove approximation properties. The global ghost penalty method is the most common approach in the literature, while the patch-wise ghost penalty method has been considered in \cite{BNV22,LZ21} and is sometimes referred to as the \emph{weak aggregation} approach.

Different realisations of the ghost penalty stabilisation method are possible. These have essentially the same properties and decompose into facet contributions:
\begin{subequations} \label{eq:sh}
\begin{align}\label{eq:sh1}
  \Sh(u,v) = \sum_{F \in \Fhgp} \gamma ~ s_{h,F}(u,v)
\end{align}
where $\gamma > 0$ is a corresponding stabilisation parameter. Two possible and popular choices for $s_{h,F}$ are:
\begin{align}\label{eq:sh2}
  s_{h,F}^{1}(u,v)= \sum_{\ell = 1}^k (h_F^{2\ell -1} \jump{ \partial_{\nF}^\ell u}, \jump{ \partial_{\nF}^\ell v})_{F}, \qquad
  s_{h,F}^{2}(u,v)= (h_F^{-2} \jump{ \Pi_{\Th} u}_{\omega_F}, \jump{\Pi_{\Th} v}_{\omega_F})_{\omega_F}.
\end{align}
\end{subequations}
Here, $\Pi_{\Th}$ denotes the element-wise $\xLtwo$ projection onto $\PP^k(\Th)$, $\omega_F = \Int( \overline{T}_1 \cup \overline{T}_2$) denotes the element aggregation to a facet $F \in \Fhgp$, $F= \partial T_1 \cap \partial T_2$.
The facet (volumetric) patch jump $\jump{u}_{\omega_F}$ of a polynomial $u_h\in\PP^k(\Th)$ is defined as
\begin{equation*}
\jump{u_h}_{\omega_F}|_{T_i} = u_h|_{T_i} - \Ex^P(u_h|_{T_j}),\quad\text{for }i,j\in\{1,2\}\text{ and } i\neq j,
\end{equation*}
where $\Ex^P$ denotes the canonical extension of a polynomial from $T$ to $\tilde \Omega$, i.e.
\begin{equation*}
    \Ex^P:\PP^k(T)\rightarrow \PP^k(\tilde\Omega)\ \text{such that}\ \Ex^Pv|_T=v \ \text{for} \ v \in \PP^k(T).
\end{equation*}
To keep the discussion simple, we only use the \emph{direct} version $s_{h,F} = s_{h,F}^2$ introduced in \cite{Pre18} in the following.

Essentially, the ghost penalty stabilisation ensures control of finite element functions on cut elements by
borrowing it from interior neighbours. For a more detailed discussion on ghost
penalties and different realisations, we refer to \cite{LO19,GM19}. The main
required property of the ghost penalty operator is that a stabilised version of
Lemma~\ref{lemma.shape-regular-norm-equivalence} holds.
\begin{lmm}\label{lemma.ghost-pebnalty-norm-equivalence}
Under Assumption~\ref{assumption.path-path}, we have for $v_h\in\PP^k(\Th)$ that
\begin{equation}\label{eq.norm-equivalence.gp}
  \nrm{\nabla v_h}_{\OT}^2 \simeq \nrm{\nabla v_h}_{\Omega}^2 + \Sh(v_h, v_h)
\end{equation}
\end{lmm}
\begin{proof}
We refer to \cite[Lemma~5.2]{LO19}.
\end{proof}

\subsection{Summary of approaches} 
\begin{figure}
  \centering
\tikzset{
  arrow/.style  = {thick,->,>=stealth},
  basic/.style  = {draw, drop shadow, font=\sffamily, thin},
  input/.style   = {basic, text width=3cm, rounded corners=2pt, align=center, fill=green!40!gray!30, rectangle},
  shapereg/.style = {basic, diamond, aspect=4, align=center, rounded corners=8pt, fill=orange!30!purple!60},
  stabil/.style = {basic, text width=4cm, align=center, fill=blue!60!green!50, rectangle},
  discr/.style = {basic, text width=4cm, rounded corners=6pt, align=center, fill=blue!60!green!70!black!20, rectangle},
  method/.style = {basic, text width=3.18cm, rectangle, rounded corners=4pt, align=center, fill=blue!70!green!20}
}

\begin{tikzpicture}[
    level 1/.style={sibling distance=5.3cm},edge from parent/.style={<-,draw,thick},>=latex]
  \tikzstyle{every node}=[font=\small]
\node (in) [input] {Input mesh $\Th$};
  
  \node (choice) [shapereg, below=.5cm of in] {Is $\Th$ cut-shape-regular?};

  \node (stabil) [stabil, right=2.8cm of choice] {Choice of stabilisation};

  \node (ag) [method, below left=1cm of stabil, xshift=.7cm] {\eqref{case:ag} $\Th=\Thag$};
  \node (wgp) [method, below right=1cm of stabil, xshift=-.7cm] {\eqref{case:pgp} $\Fhgp=\Fhgpmin$};
  \node (gp) [method] at ($(ag)!0.5!(wgp)$) {\eqref{case:gp} $\Fhgp=\Fhgpstar$};
  
  \node (methodgp) [discr, below=.8cm of gp] {Stabilised discretisation\\$\Thag\neq\Th$, $\Sh(\cdot,\cdot)\neq 0$.};
  \node (methodshapereg) [discr] at (methodgp -| choice) {Unstabilised discretisation \\$\Thag=\Th$, $\Sh(\cdot,\cdot)\equiv 0$.};

  \draw [arrow] (in) -- (choice);
  \draw [arrow] (choice) -- node[anchor=south] {no} (stabil);
  \draw [arrow] (choice) -- node[anchor=east] {yes} (methodshapereg);
\draw [arrow] (stabil) -- (ag);
  \draw [arrow] (stabil) -- (gp);
  \draw [arrow] (stabil) -- (wgp);

  \draw [arrow] (ag) -- (methodshapereg);
  \draw [arrow] (gp) -- (methodgp);
  \draw [arrow] (wgp) -- (methodgp);
\end{tikzpicture}  
   \caption{Stabilisation choices dependent on the input mesh and resulting mesh and facet sets.}
  \label{fig.flow-diagramm}
\end{figure}
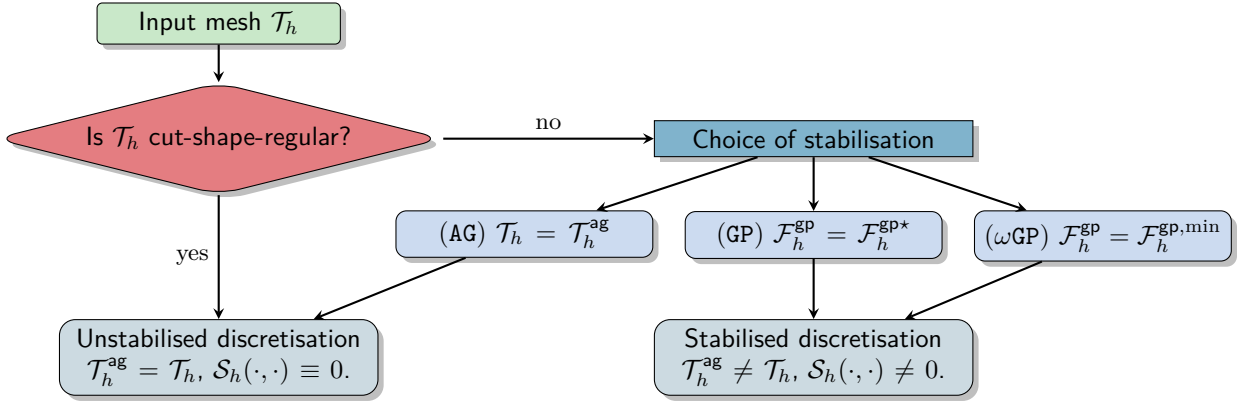

\noindent
We finalise this section by collecting the three different methods of stabilisation under consideration:
\begin{align}
  \text{ \tag{\texttt{GP}} \label{case:gp}
    \begin{minipage}{0.88\textwidth} 
        The ghost penalty method uses a stabilisation term, given in \eqref{eq:sh}, acting on $\Fhgp=\Fhgpstar$.
        The element clusters are then used solely in the analysis for interpolation. 
        In this case, the shape regularity of the aggregated elements (and the nestedness of finite element spaces) can be exploited to obtain optimal approximation results.
    \end{minipage}
  } \\[0.4ex]
  \text{ \tag{$\omega$\texttt{GP}} \label{case:pgp}
      \begin{minipage}{0.88\textwidth} 
          For the patch-wise ghost penalty method, the same stabilisation term \eqref{eq:sh} is used over the minimal set of faces. The element clusters are used to reduce the regions where stabilisation is applied, as stability for bad-cut elements can be supported from within each cluster by choosing $\Fhgp=\Fhgpmin$.
    \end{minipage}
  } \\[0.4ex]
  \text{ \tag{\texttt{AG}} \label{case:ag}
    \begin{minipage}{0.88\textwidth} 
      For the element aggregation method, mesh elements are merged, resulting in an active mesh $\Th  = \Thag$ that is \emph{cut-shape-regular} in the sense of Definition~\Rref{def.shaperegcut}. Hence, no ghost penalty-like stabilisation is needed, and we set $\Sh(\cdot,\cdot)\equiv 0$. 
  \end{minipage}
  } 
\end{align}
We note that one could also think about a situation where one directly
starts with a \emph{cut-shape-regular} mesh $\Th$ and then considers an
unstabilised discretisation. In the analysis below, we will take this
viewpoint for the case (\texttt{ag}), i.e., we simply assume to be given a
\emph{cut-shape-regular} mesh $\Th$ for which no stabilisation is required.
Given a particular mesh, we illustrate the choices available in Figure~\ref{fig.flow-diagramm} and the resulting consequences for our mesh notation.

\section{Error Analysis} \label{sec:errana}
In this section, we analyse the unfitted DG and Trefftz DG methods under the assumption of exact geometry handling.
The analysis for the unfitted DG method has already been treated in \cite{GM19}. We repeat the analysis with slight generalisations concerning mesh assumptions for inverse inequalities, the aggregated DG formulation and a patch-wise ghost penalty formulation. In particular, we generalise the approximation result by a special interpolation operator in Section~\ref{sssec:approx}. The particular interpolation operator allows us to conveniently extend the analysis further for the unfitted Trefftz DG method.

\subsection{Error Analysis of the Unfitted DG and Trefftz DG Methods with Exact Geometry}
\label{sec.error-exact-geometry}
In this section, we will analyse the unfitted DG as well as the unfitted Trefftz
DG method introduced above under the assumption that we are given a smooth
particular solution $u_p\in \xHn{\ell}(\Omega),~\ell\geq 2$ of \eqref{eqn.strong-poisson.bulk}. We homogenise the problem
a priori, i.e.\@{} instead of $u$, the solution to \eqref{eqn.strong-poisson}, we
look for $u$, the solution to the following variant:
Find $u\colon \Omega \to \mathbb{R}$ such that
\begin{subequations}\label{eqn.strong-poisson.hom}
  \begin{align}
    -\Delta u &= 0 \quad \textnormal{in } \Omega \tag{\ref{eqn.strong-poisson}a\textsuperscript{hom}}\\
    u &= \ghom = g-u_p \quad \textnormal{on } \Gamma. \tag{\ref{eqn.strong-poisson}b\textsuperscript{hom}}
  \end{align}
\end{subequations}
Hence, no additional homogenisation is needed in the numerical methods. The
corresponding DG and Trefftz DG method (both with right-hand side $\Lh(0, \ghom; v_h)$) 
will be denoted as \eqrefh{eqn.poisson-poblem.dg} and \eqrefh{eqn.poisson-poblem.trefftz}. A numerical solution for the more generic case that no suitable $u_p\in \xHn{\ell}(\Omega)$ is known will be discussed in Section~\ref{sec.further-methods}.

In the remainder of this work, we make the following assumptions.
Firstly, we assume that the computational mesh $\Th$ fulfils
Assumption~\Rref{assumption.mesh2}. Second, either the stabilisation bilinear
form $S_h(\cdot,\cdot)$ is chosen as described in
Section~\ref{sec.stabilization:subset.ghost-penalty} and
Assumption~\Rref{assumption.path-path} is valid, or
$\Th$ is already \emph{cut-shape-regular} so that $S_h(\cdot,\cdot)=0$. In the
latter case, we set $\Thag = \Th$, i.e.\@{} all patches considered below are trivial
patches only consisting of one element in $\Th$. 

\subsubsection{Trace Inequalities and Coercivity}
We gather several trace estimates which hold due to Assumption~\Rref{assumption.mesh2}. 
\begin{lmm}[Continuous trace inequality]\label{lem:traceineq1}
Let $T \in \Th$. For all $v \in \xHn{1}(T)$, we have
\begin{equation}\label{eqn.traceineq1}
  \Vert v \Vert_{\partial T}^2 \lesssim h_T^{-1} \Vert v \Vert_{T}^2 + h_T \Vert \nabla v \Vert_{T}^2.
\end{equation}
\end{lmm}
\begin{proof}
  For meshes that fulfil Assumption~\Rref{assumption.mesh1}, the proof can be
  found in \cite[Lemma~4.7]{CDG21}. Let $T \in \Th$ with subelement $\{T'\}$ as
  in Assumption~\Rref{assumption.mesh2}. We can then divide the boundary of $T$
  into the subelement contribution and apply the trace inequality for the
  subelement $T'$,
  \begin{equation*}
     \Vert v \Vert_{\partial T}^2
     = \sum_{T'} \Vert v \Vert_{\partial T' \cap \partial T}^2
     \lesssim \sum_{T'} h_{T'}^{-1} \Vert v \Vert_{T'}^2 + h_{T'} \Vert \nabla v \Vert_{T'}^2
     \lesssim h_T^{-1} \Vert v \Vert_{T}^2 + h_T \Vert \nabla v \Vert_{T}^2.
     \vspace*{-1.1cm}
  \end{equation*} 
  \vspace*{0.1cm}
\end{proof}
\begin{lmm}[Discrete trace inequality]\label{lemma.discrete-trace}
  Let $T \in \Th$. For all $v_h\in\PP^k(T)$, we have
  \begin{equation}
    \nrm{v_h}_{\partial T}^2 \leq c k^2 h_T^{-1} \nrm{v_h}_T^2,
  \end{equation}
  with a constant $c>0$, independent of the local mesh size $h_T$ and the order $k$.
\end{lmm}
\begin{proof}
  The claim follows from Lemma~\ref{lem:traceineq1} and an inverse inequality
  applied on the subelements $T'$ of $T$, cf.\@{} Assumption~\Rref{assumption.mesh2},
  \begin{equation}\label{eqn.inverse.estimate}
    \Vert \nabla v_h \Vert_{T}^2 = \sum_{T'} \Vert \nabla v_h \Vert_{T'}^2 \leq
      c k^4 h_T^{-2} \sum_{T'} \Vert v_h \Vert_{T'}^2
    = k^4 h_T^{-2} \Vert v_h \Vert_{T}^2  
  \end{equation}
  for all $v_h \in \PP^k(T)$, and a constant $c>0$ independent of $h$ and $k$.
  See~\cite[Corollary~4.24]{CDG21} for details on the inverse inequality
  applied on $T'$ under our mesh assumptions. 
\end{proof}

\begin{lmm}[Continuous unfitted trace inequality]\label{lem:traceineq2}
  Let $T \in \Th$ and $T^\Gamma = T \cap \Gamma$ with $\meas_{d-1}T^\Gamma>0$. Then, there holds for any $v \in \xHn{1}(T)$:
  \begin{equation}\label{eqn.traceineq2}
    \Vert v \Vert_{T^\Gamma}^2 \lesssim h_T^{-1} \Vert v \Vert_T^2 + h_T \Vert \nabla v \Vert_{T}^2.
  \end{equation}
\end{lmm}
\begin{proof}
  Similar proofs to similar statements are given e.g. in \cite{HH02,DER14,BE18}. For the framework at hand, the proof of \cite[Lemma 5.2]{LX21} can be used with only slight adaptations in the setting. For completeness, we included the proof in Lemma~\ref{lem:traceineq2:proof}.
\end{proof}

In conjunction with the problem \eqref{eqn.poisson-poblem.dg}, we introduce the following (semi-)norms for the analysis below
\begin{subequations}
\begin{align}
  \tnrmA{v}^2 &\coloneqq \nrm{\nabla v}_{\Omega}^2 + \nrm{h^{-1/2}[v]}_{\Fh\cap\Omega}^2 + \nrm{h^{1/2}\mean{\dnF v}}_{\Fh\cap\Omega}^2 + \nrm{h^{-1/2}v}_{\Gamma}^2 + \nrm{h^{1/2}\dn v}_{\Gamma}^2,\\
  \snrmS{v}^2 &\coloneqq \Sh(v,v),\\
  \tnrmB{v}^2 &\coloneqq \tnrmA{v}^2 + \snrmS{v}^2.
\end{align}
\end{subequations}

\begin{lmm}\label{lemma.poisson.dg}
$\Bh(\cdot,\cdot)$ is coercive and $\Ah(\cdot,\cdot)$, $\Bh(\cdot,\cdot)$ and $\Lh(\cdot)$ are continuous, i.e., for all $u\in \xHtwo(\OT)\oplus \PP^k(\Th)$ and $v_h\in \PP^k(\Th)$ there holds
   \begin{align}
       \qquad&&&& \Ah(u, v_h) &\lesssim \tnrmA{u}\tnrmA{v_h}, &&& \Lh(v_h)&\lesssim\tnrmA{v_h} &&&&\qquad \label{eq.lemma.poisson.dg.1} \\
       &&&& \Bh(v_h, v_h) &\gtrsim \tnrmB{v_h}^2 & \text{ and } && \Bh(u, v_h) &\lesssim \tnrmB{u}\tnrmB{v_h}, \label{eq.lemma.poisson.dg.2} 
   \end{align}
   if $\beta>0$ is sufficiently large, independent of $h$ and $k$. Finally, the problem \eqref{eqn.poisson-poblem.dg} admits a unique solution.
\end{lmm}

\begin{proof}
The ghost-penalty stabilised case is covered by \cite[Proposition~2.6]{GM19} for
shape-regular (but not necessarily cut-shape-regular) simplicial meshes. We briefly repeat the arguments for this setting. The continuity estimates follow immediately from the Cauchy-Schwarz inequality. For the coercivity, we have with the Cauchy-Schwarz and weighted Young's inequality
\begin{align*}
  \Ah(v_h, v_h) &= \nrm{\nabla v_h}_{\Omega}^2 - 2(\dn v_h, v_h)_\Gamma + \beta  \nrm{ h^{-\frac12} v_h}_\Gamma^2 
    - 2(\mean{\dnF v_h},\jump{v_h})_{\Fh\cap\Omega} + \beta \nrm{ h^{-\frac12}\jump{v_h}}_{\Fh\cap\Omega}^2\\
    &\geq \nrm{\nabla v_h}_{\Omega}^2 - \epsilon\nrm{h^{\frac{1}{2}}\dn v_h}_\Gamma^2 + \Big(\beta\!-\!\frac{1}{\epsilon}\Big) \nrm{h^{-\frac{1}{2}}v_h}_\Gamma^2
      - \epsilon\nrm{h^{\frac{1}{2}}\mean{\dnF v_h}}_{\Fh\cap\Omega}^2 + \Big(\beta\!-\!\frac{1}{\epsilon}\Big) \nrm{h^{-\frac{1}{2}}\jump{v_h}}_{\Fh\cap\Omega}^2.
\end{align*}
We then observe with $\sigma_i=\dnF|_{T_i} v_h$ and the discrete trace inequality in Lemma~\ref{lemma.discrete-trace} that
\begin{equation}\label{eqn.bound-mean-normal-deriv}
  \nrm{h^{\frac{1}{2}}\mean{\dnF v_h}}_{\Fh\cap\Omega}^2
    \leq \sum_{F\in\Fh} h_F (\nrm{\sigma_1}_F^2 + \nrm{\sigma_2}_F^2)
    \sum_{T\in\Th} h_T \nrm{\sigma}_{\partial T}^2
    \leq \sum_{T\in\Th} c_1 k^2 \nrm{\nabla v_h}_{T}^2
    \leq c_1 k^2 \nrm{\nabla v_h}_{\OT},
\end{equation}
with a constant $c_1>0$ independent of the local mesh size $h_T$ and order $k$.
With a cut version of the discrete trace inequality, obtained by combining Lemma~\ref{lem:traceineq2} and \eqref{eqn.inverse.estimate}, we also have
\begin{equation}\label{eqn.bound-bnd-normal-deriv}
  \nrm{h^{\frac{1}{2}}\dn v_h}_\Gamma \leq c_2 k^2 \nrm{\nabla v_h}_{\OT},
\end{equation}
with a constant $c_2>0$ independent of the mesh size (field) $h$ and order $k$.
In the cut-shape-regular case, we can use Lemma~\ref{lemma.shape-regular-norm-equivalence} to bound the right-hand side of \eqref{eqn.bound-mean-normal-deriv} and \eqref{eqn.bound-bnd-normal-deriv} by a norm on $\Omega$, and in the ghost penalty stabilised case, we can use Lemma~\ref{lemma.ghost-pebnalty-norm-equivalence} to the same effect. Consequently, we have
\begin{equation*}
  \nrm{h^{\frac{1}{2}}\mean{\dnF v_h}}_{\Fh\cap\Omega}^2 + \nrm{h^{\frac{1}{2}}\dn v_h}_\Gamma^{2}
  \leq c k^2 (\nrm{\nabla v_h}_{\Omega}^2 + \Sh(v_h, v_h)),
\end{equation*}
where we recall that in the cut shape regular case the ghost penalty terms are
empty, i.e.\@{} $\Sh(\cdot,\cdot) = \snrmS{\cdot}^2 \equiv 0 $.
Combining this with the above estimate, choosing $\epsilon\leq\frac{1}{2ck^2}$ and subsequently $\beta \geq 4ck^2$, we have
\begin{align*}
  \Bh(v_h, v_h) \geq (1 - \epsilon c k^2)(\nrm{\nabla v_h}_\Omega^2 + \snrmS{v_h}^2)
    + (\beta - \epsilon^{-1})\Big( \nrm{h^{-\frac{1}{2}}v_h}_\Gamma^2 + \nrm{h^{-\frac{1}{2}}\jump{v_h}}_{\Fh\cap\Omega}^2  \Big) \geq \frac{1}{2}\tnrmB{v_h}^2.
\end{align*}
\end{proof}

In the following we will always assume that $\beta$ is sufficiently large so that Lemma~\ref{lemma.poisson.dg} holds.

\begin{crllr} \label{corollary.poisson.tdg}
  The problem in \eqrefh{eqn.poisson-poblem.trefftz} admits a unique solution and \eqref{eq.lemma.poisson.dg.1}--\eqref{eq.lemma.poisson.dg.2} also holds for $v_h \in \TT^k(\Th)$.
\end{crllr}
\begin{proof}
  Coercivity and continuity as stated in Lemma~\ref{lemma.poisson.dg} are directly inherited on the subspace $\TT^k(\Th) \subset \PP^k(\Th)$, which follows per definition of the DG and Trefftz spaces in \eqref{eqn.dg-space} and \eqref{eqn.trefftz-space}.
\end{proof}

\subsubsection{Céa-type Quasi-Best Approximation Results}

\begin{lmm} \label{lem.poisson.dg+tdg}
 Let $u\in \xHn{1}(\Omega) \cap \xHtwo(\Th)$ be the solution to
 \eqrefh{eqn.strong-poisson} with $\ghom \in \xHn{\frac{1}{2}}(\Gamma)$.
 Furthermore, let either $V_h=\PP^k(\Th)$ and $u_h\in V_h$ be the solution to
 \eqrefh{eqn.poisson-poblem.dg}, or let
 $V_h=\TT^k(\Th)$ and $u_h\in V_h$ be the solution to
 \eqrefh{eqn.poisson-poblem.trefftz}. Then
 \begin{equation}
   \tnrmA{u - u_h} + \snrmS{u_h}  \lesssim \inf_{v_h \in V_h} \tnrmA{u - v_h} + \snrmS{v_h}. 
 \end{equation}
\end{lmm}

\begin{proof}The proof for the DG case has essentially been given in \cite[Theorem 2.10]{GM19}.
With the triangle inequality we have for arbitrary $v_h \in V_h$ that
\begin{equation*}
  \tnrmA{u - u_h} + \snrmS{u_h} \leq\tnrmA{u - v_h} + \tnrmA{v_h - u_h} + \snrmS{u_h - v_h} + \snrmS{v_h} \!\lesssim\! \tnrmA{u - v_h} + \tnrmB{v_h - u_h} + \snrmS{v_h}.
\end{equation*}
With $w_h \coloneqq u_h - v_h$, Lemma~\ref{lemma.poisson.dg} and Corollary~\ref{corollary.poisson.tdg} and $\Bh(u_h,v_h)=\Ah(u,v_h)~\forall v_h \in V_h$, we see
\begin{align*}
  \tnrmB{w_h}^2 &\lesssim \Bh(u_h - v_h, w_h) = \Ah(u - v_h, w_h) - \Sh(v_h, w_h)\\
    &\lesssim \tnrmA{ u - v_h}\tnrmA{w_h} + \snrmS{v_h}\snrmS{w_h}
    \lesssim (\tnrmA{ u - v_h} + \snrmS{v_h})(\tnrmA{w_h} + \snrmS{w_h})
\end{align*}
Dividing by $\tnrmA{w_h} + \snrmS{w_h} \simeq \tnrmB{w_h}$ and combining with the previous triangle inequality concludes the proof.
Note that all steps are valid for the unfitted DG as well as for the unfitted Trefftz DG case.
\end{proof}

\subsubsection{Approximation} \label{sssec:approx}

For the approximation results, we will use averaged Taylor polynomials on adapted domains as the interpolation operator, cf.\@{} \cite[Section 4.1]{BS08}. We repeat its crucial properties and formulate them in a suitable setting for the following analysis. 

\begin{lmm} \label{lem.avg.Taylor}
Let $T$ be a domain with Lipschitz boundary and assume that there are two balls $b_T$ and $B_T$ with $b_T \subset T \subset B_T$ and associated diameters $h_b$ and $h_B$ so that $h_B/h_b = \sigma_T \lesssim 1$. We define the interpolation operator $\pi_T^m: \xLone(b_T) \to \mathcal{P}^m(B_T)$ as the operator realizing the averaged Taylor polynomial of degree $m$ by averaging over $b_T$, see~\cite{BS08}.
\begin{subequations}
\begin{enumerate}
    \item For $v\in \xWn{{m+1,p}}(B_T)$ and $p\geq1$ it holds that
\begin{equation} \label{eqn.avgtaylor.approx}
    \snrm{v - \pi_T^m v}_{\xWn{{n,p}}(T)}\leq \snrm{v - \pi_T^m v}_{\xWn{{n,p}}(B_T)} \lesssim h_B^{m+1-n}\snrm{v}_{\xWn{{m+1,p}}(B_T)}\qquad\text{for } n=1,\dots,m+1,
\end{equation}
with a constant only depending on $\sigma_B$, $m$ and $n$.
\item  For $v \in \xWn{{|\alpha|,1}}(T)$ and $\alpha \in \mathbb{N}^d$ such that $|\alpha| \leq m-1$ there holds
\begin{equation}\label{eqn.avgtaylor.commutes}
    D^\alpha \pi_T^m v = \pi_T^{m-|\alpha|} D^\alpha v.
\end{equation}
\end{enumerate}
\end{subequations}
\end{lmm}
\begin{proof}
See \cite[Section 4.1]{BS08}, especially \cite[Lemma 4.3.8 \& Proposition 4.1.17]{BS08}, since $B_T$ is star-shaped.
\end{proof}
The second property of commuting interpolation and differentiation is crucial in the context of the Trefftz DG subspace. Note that in \eqref{eqn.avgtaylor.commutes} $v$ and $D^\alpha v$ are effectively only evaluated on the small ball $b_T$.
It specifically implies that for $u \in \xHn{m}(T)$ with $\Delta u = 0$ (on $b_T$) and $m\geq 2$ there holds $\Delta \pi_T^m u = \pi_T^{m-2} \Delta u = 0$ on each element (for $m < 2$ there trivially holds $\Delta \pi_T^m u = 0$). In other words, the averaged Taylor polynomial of a harmonic function is also harmonic.

As a consequence of Lemma~\ref{lmm.ass1ass2} and the fact that the aggregated elements in
$\Thag$ again fulfil Assumption~\ref{assumption.mesh2} we make the following
observation in preparation of the subsequent lemma.
\begin{crllr} \label{crllr.avgtaylor}
  Let $T \in \Thag$. For $b_T$ the largest ball in $T \cap
  \Omega$ and $B_T$ the smallest ball with $B_T \supset T$ \eqref{eqn.avgtaylor.approx} and \eqref{eqn.avgtaylor.commutes} hold.
\end{crllr}
\begin{lmm}\label{lem.approx.TTh}
  For $u \in \xHn{m}(\Omega)$ with $\Delta u = 0$ and $l = \min \{m-1, k\}$ there holds
  \begin{align}
\inf_{v_h \in \TT^k(\Th)} \left( \tnrmA{ u - v_h} \!\!+ \snrmS{v_h} \! \right)\!\leq\!  \inf_{v_h \in \TT^k(\Thag)} \left( \tnrmA{ u - v_h} \!\! + \snrmS{v_h} \! \right) \nonumber \!\lesssim\! 
\left( { \textstyle\sum_{T \in \Th} } h_T^{2l} \Vert \mathcal{E} u \Vert_{\xHn{{l+1}}(T)}^2 \right)^{\!1/2}  \!\!\!\!
\lesssim 
h^{l} \Vert u \Vert_{\xHn{{l+1}}(\Omega)},
  \end{align}
  where the constant depends on the maximum number of subelements $T'$ in each
  element, $m_T$, cf.\@{} Assumption~\ref{assumption.mesh2}, and the maximum number of elements in a patch $n_{max}$, cf., Lemma~\ref{lemma.patch-size}.
\end{lmm}
\begin{proof}
  We first recall that in the case where $\Th$ is \emph{cut-shape-regular}
  we have $\Thag = \Th$.
  
  The first inequality is obvious due to $\TT^k(\Thag) \subset \TT^k(\Th)$.
For the second, we will find an interpolator $v_h \in \TT^k(\Thag)$ with the
  corresponding bound. Before we discuss the construction in more detail we want
  to stress that the interpolator will be constructed w.r.t. the finite element
  space on $\Thag$ which is smaller than the one used in the discretisation
  \emph{only} in the case where $\Th$ is not \emph{cut-shape-regular}
  (otherwise $\Thag=\Th$ holds).
  
  First, we bound all different norm contributions by $\xHn{j}$(T)-semi-norms on elements $T$ in the active mesh $\Th$.
To this end, we recall the trace inequalities \eqref{eqn.traceineq1} and \eqref{eqn.traceineq2}.
Applying these estimates to all norm contributions of $\tnrmA{\cdot}$ to $v \in \xHn{m}(\Th)$ yields
\begin{align*}
  \tnrmA{v}^2 
  &\lesssim \sum_{T \in \Th} \nrm{\nabla v}_{T \cap \Omega}^2 
    + h_T^{-1} \nrm{v}_{\partial T}^2 
    + h_T \nrm{\nabla v}_{\partial T}^2
    + h_T^{-1} \nrm{v}_{\Gamma \cap T}^2 
    + h_T \nrm{\nabla v}_{\Gamma \cap T}^2 \\
  &\lesssim \sum_{T \in \Th} \snrm{v}_{\xHn{1}(T)}^2 
    + h_T^{-2} \nrm{v}_T^2 + h_T^2 \snrm{v}_{\xHtwo(T)}^2
  \lesssim \sum_{T \in \Th} \sum_{j = 0,1,2} h_T^{2(j-1)} \snrm{v}_{\xHn{j}(T)}^2.
\end{align*}
For the construction of the interpolant $v_h = \inter^k u$, we apply the
averaged Taylor polynomial from Lemma~\ref{lem.avg.Taylor} patch-wise for all $\omega
\in \Thag$  so that $\inter^k u|_{\omega} = \pi_{\omega}^k u$. Making use of
Corollary~\ref{crllr.avgtaylor} the averaged Taylor polynomial $\inter^k u|_{\omega}$
is based only on $u|_{b_\omega}$ with a ball $b_\omega\subset \omega\cap\Omega$. By construction $v_h$ depends only on $u$ in $\Omega$ and hence with \eqref{eqn.avgtaylor.commutes} is harmonic, i.e.\@{} $v_h\in\TT^k(\Thag) \subset \TT^k(\Th)$. 

We cannot directly plug in the approximation error bound for $u$, as $u$ is not defined on $\OT$. We hence make use of a linear extension operator $\Ex\colon \xHn{m}(\Omega)\rightarrow \xHn{m}(\OT)$, for $m\geq 0$, for which it holds
\begin{equation*}
    \Ex u|_{\Omega} = u \text{ and } \nrm{\Ex u}_{\xHn{m}(\OT)} \lesssim \nrm{u}_{\xHn{m}(\Omega)},
\end{equation*}
see for example \cite[Section VI.3]{Ste70}. For the interpolant of the extension
$\Ex u$ on the entire aggregated element $\omega$, we choose the same averaged Taylor polynomial, so that $\pi_h^k\Ex u = \pi_h^k u \in \TT^k(\Thag)$. We note that we constructed the averaged Taylor polynomial especially so that $\inter^k$ only depend on values in $\Omega$.
We set $v_h = \inter^k (\Ex u)$. Applying \eqref{eqn.avgtaylor.approx} gives for $v \in \xHn{{m+1}}(\Omega)$ with $\Delta v = 0$ and each $\omega \in \Thag$ the estimate
$$\vert \Ex v - \pi_h^k\Ex v \vert_{\xHn{j}(\omega)} \lesssim h^{m+1-j}_B \vert \Ex v \vert_{\xHn{{m+1}}(B_\omega)}$$
for $ j,m \leq k$. Hence, we have (with a finite overlap argument for the domains $B_\omega, \omega \in \Thag$)
\begin{multline*}
  \tnrmA{u - v_h}^2  \leq \tnrmA{\Ex u - \inter^k \Ex u }^2
  \lesssim \sum_{\omega \in \Thag} \sum_{j = 0,1,2} h_\omega^{2(j-1)} \left( \snrm{\Ex u - \inter^k \Ex u  }_{\xHn{j}(\omega)}^2 \right)\\
  \lesssim \sum_{\omega \in \Thag} h_\omega^{2l} \left( \snrm{\Ex u}_{\xHn{{l+1}}(B_\omega)}^2 \right) 
  = h_\omega^{2l} \snrm{\Ex u}_{\xHn{{l+1}}(\OT)}^2  \left( \lesssim h^{2l} \snrm{u}_{\xHn{{l+1}}(\Omega)}^2 \right)
\end{multline*}
For the ghost-penalty semi-norm, let us consider two aligned elements $T_1,\ T_2$ and an aggregated element from a single facet patch
$\omega_F=\Int(\overline{T}_1\cup \overline{T}_2)$ and $\Tw = \{T_1,T_2\}$. 
We denote with $\Pi_{\omega_F}$ the $\xLtwo$ projection onto the polynomial space over the domain $\omega_F$, and by $\Pi_{\Tw}$ the element-wise projection onto the broken polynomial space over the elements $\Tw$.
Let us denote $v_{h,i} = \Ex^P \Pi_{T_i}v_h  
\in \PP^k(\omega_F)$, $i=1,2$.
Then for any $v \in
\PP^k(\omega_F)$ 
\begin{align*}
    h_F^2 \vert{v_h}\vert_{\Sh, \omega_F}^2
= \nrm{v_{h,1} - v_{h,2}}_{\omega_F}^2
        &\leq 2(\nrm{v_{h,1} - v}_{\omega_F}^2 + \nrm{v - v_{h,2}}_{\omega_F}^2)\\
        &= 2 (\nrm{v_{h,1} - v}_{T_1}^2 + \nrm{v_{h,1} - v}_{T_2}^2 + \nrm{v - v_{h,2}}_{T_1}^2 + \nrm{v - v_{h,2}}_{T_2}^2)\\
        &\lesssim \nrm{\Pi_{T_1} v_{h} - v}_{T_1}^2 + \nrm{v - \Pi_{T_2}  v_{h}}_{T_2}^2
        = \nrm{ \Pi_{\Tw} v_{h} - v}_{\omega_F}^2
\end{align*}
where the second last step follows using the shape regularity of the mesh, by which we
can bound the norm of the discrete function on $T_1$ by the norm on $T_2$ and
vice versa as the arguments are polynomials (not only element-wise) on the aggregated element. 
Let $v= \Pi_{\omega_F}(\Ex u)$ and recall $v_h=\pi^k_h(\Ex u)$ then
\begin{align*}
  \Pi_{\Tw} v_h - v = (\inter^k - \Pi_{\omega_F}) \Ex u  =  (\inter^k - \id) \Ex u + (\id - \Pi_{\omega_F}) \Ex u.
\end{align*}
The operators $\inter^k$, $\Pi_{\omega_F}$
have the usual optimal
approximation bounds.
We hence obtain
\begin{equation} \label{eqn.ghost-penalty-consitency}
  \!\!\vert{v_h}\vert_{\Sh,\omega_F} \! \lesssim h_F^{-1}\nrm{ \Pi_{\Th} v_h - v}_{\omega_F}
  \! \lesssim h_F^l \nrm{\Ex u}_{\xHn{{l+1}}(\omega_F)} 
   ~\Rightarrow \vert{v_h}\vert_{\Sh}
    \lesssim \left( {\textstyle \sum_{F \in \Fhgp} } h_F^{2l} \nrm{\Ex u}_{\xHn{{l+1}}(\omega_F)}^2 \right)^{\frac12}
   \!\!\lesssim h^l \nrm{u}_{\xHn{{l+1}}(\Omega)}.\!
\end{equation}
We note that in the case of the patch-wise ghost-penalty operator, i.e., $\Fhgp= \Fhgpmin$, we have by construction $\snrm{\pi_h^m \Ex u}_{\Sh} = 0$ for $m\leq k$  as $\pi_h^m$ maps onto $\mathbb{P}^m(\Thag)$ which is in the kernel of $| \cdot |_{\Sh}$ in the case of the patch-wise ghost penalty.
\end{proof}

\begin{crllr}\label{cor:trefftz:approx}
  For $u \in \xHn{m}(\Omega)$ with $\Delta u = 0$ and $l = \min \{m-1, k\}$ there holds
  \begin{equation}
    \inf_{v_h \in \PP^k(\Th)} \left( \tnrmA{u - v_h} + \snrmS{v_h} \right) \lesssim h^{l} \Vert u \Vert_{\xHn{{l+1}}(\Omega)}.
  \end{equation}
\end{crllr}
\begin{proof}
  Follows from $ \TT^k(\Thag) \subset \TT^k(\Th) \subset \PP^k(\Th)$ and Lemma~\ref{lem.approx.TTh}.
\end{proof}

\subsubsection{A Priori Error Bounds}

\begin{crllr}\label{corrolary.energy-estimate}
  Let $u \in \xHn{m}(\Omega)$ be the solution to \eqrefh{eqn.strong-poisson} with $\ghom\in \xHn{\frac{1}{2}}(\Gamma)$ and $u_h \in V_h$ be the solution to
  \eqrefh{eqn.poisson-poblem.dg} or \eqrefh{eqn.poisson-poblem.trefftz} with $V_h=\PP^k(\Th)$ or $V_h=\TT^k(\Th)$, respectively. Then, there holds for $l = \min\{m-1,k\}$
  \begin{equation}
    \tnrmA{u - u_h} + \snrmSg{u_h} \lesssim h^{l} \Vert u \Vert_{\xHn{{l+1}}(\Omega)}.
  \end{equation}
\end{crllr}
\begin{proof}
  Follows from Lemma~\ref{lem.poisson.dg+tdg} and Lemma~\ref{lem.approx.TTh}.
\end{proof}

\begin{thrm} \label{thm:l2est}
  We assume $\Omega$ to be sufficiently smooth or convex such that
  $\xLtwo(\Omega)$-$\xHtwo(\Omega)$-regularity holds\footnote{This means that for $f\in L^2(\Omega)$ the solution $w$ to $-\Delta w\!=\!f$ in $\Omega$, $w\!=\!0$ on $\Gamma$ is $H^2$-regular, $w \in H^2(\Omega)$ and $\Vert w \Vert_{H^2(\Omega)}\!\lesssim\! \Vert f \Vert_{L^2(\Omega)}$.} and assume that $\Th$ is quasi-uniform, i.e. $h_T \simeq h~\forall T \in \Th$. Furthermore, let $u \in
  \xHn{m}(\Omega),~m\geq 2$ be the solution to \eqrefh{eqn.strong-poisson} with $\ghom\in \xHn{\frac{1}{2}}(\Gamma)$ and $u_h \in V_h$ be the solution to
  \eqrefh{eqn.poisson-poblem.dg} or \eqrefh{eqn.poisson-poblem.trefftz}. Then, there holds for $l = \min\{m-1,k\}$ that
  \begin{equation}
    \nrm{u - u_h}_{\Omega} \lesssim h^{l+1}\nrm{u}_{\xHn{{l+1}}(\Omega)}.
  \end{equation}
\end{thrm}

\begin{proof}
Due to the elliptic regularity assumption, we have that for the auxiliary problem
\begin{alignat*}{2}
  -\Delta z &= u - u_h\quad &&\textnormal{in } \Omega, \\
          z &= 0            &&\textnormal{on } \Gamma,
\end{alignat*}
that $z\in \xHtwo(\Omega)\cap \xHn{1}_0(\Omega)$ and that $\nrm{z}_{\xHtwo(\Omega)}
\lesssim \nrm{u - u_h}_{\Omega}$. From $z\in \xHtwo(\Omega)$, it follows that
$\jump{\nabla z}\cdot\nF = 0$ and $\jump{z}=0$ on all $F\in\Fh\cap\Omega$. Due
to the symmetry and consistency of the symmetric interior penalty form
$\Ah(\cdot,\cdot)$, we have that
\begin{equation*}
  \Ah(u- u_h, z) = \Ah(z, u- u_h) = \int_{\Omega}(-\Delta z)(u-u_h)\dif \xb = \nrm{u - u_h}^2_{\Omega}.
\end{equation*}
Furthermore, we have the perturbed Galerkin-orthogonality
\begin{equation}\label{eqn.perturbed-Galerin-orthogonality}
  \Ah(u - u_h, v_h) = \Sh(u_h, v_h)\qquad\text{for all }v_h\in V_h.
\end{equation}
Now let $\pi_h^1$ be the interpolation operator by average Taylor polynomials onto $\PP^1(\Thag)=\TT^1(\Thag)\subset V_h$. We note that all piecewise linear functions are naturally harmonic. Then by \eqref{eqn.perturbed-Galerin-orthogonality} and Lemma~\ref{lemma.poisson.dg}, we have 
\begin{align*}
  \nrm{u - u_h}^2_{\Omega} = \Ah(u- u_h, z) 
  &= \Ah(u- u_h, z - \pi_h^1 z) + \Sh(u_h , \pi_h^1 z)\\
  &\lesssim  \tnrmA{u - u_h}\tnrmA{z - \pi_h^1 z} + \snrmS{u_h}\snrmS{\pi_h^1 z}\\
  &\lesssim (\tnrmA{u - u_h} + \snrmS{u_h})h \nrm{z}_{\xHtwo}
  \lesssim h(\tnrmA{u - u_h} + \snrmS{u_h})\nrm{u - u_h}_{\Omega},
\end{align*}
where the penultimate estimate follows from Lemma~\ref{lem.avg.Taylor} and the consistency of the ghost-penalty semi-norm shown in Lemma~\ref{lem.approx.TTh}. The claim then follows from Corollary~\ref{corrolary.energy-estimate}.
\end{proof}

\subsection{Unfitted DG and Trefftz DG Methods With Geometry Approximation}
\label{sec.anaysis-with-geom-error}
Unfitted finite element discretisations pose the additional challenge of geometry approximation, i.e., the demand for an accurate representation of a geometry not described through the computational mesh and the need for robust and accurate numerical integration over cut elements. Several techniques to achieve these are known in the literature; see for example, \cite{FOSS17,HSK17,Leh16,MKO13,OS16,Say15}. In our numerical examples below, we will consider an approach based on piecewise linear reference configuration and a (small) local mesh deformation ~\cite{Leh16}.
In this section, we introduce the methods with respect to a discrete
approximated geometry and carry out an error analysis based on a Strang-type
lemma. Similar techniques have been applied in \cite{LR17,Leh17} and the works by Deckelnick, Elliott, Ranner, e.g. \cite{ER12,DER14}. For ease of presentation we restrict to the case of quasi-uniform meshes, i.e. $h_T \simeq h$ for all $T\in\Th$.

\subsubsection{Geometry Approximation}
In the remainder we assume that a geometry approximation $\Omega_h$ of higher order is given, i.e.,
\begin{equation*}
  \dist(\Omega, \Omega_h) \lesssim h^{q+1},
\end{equation*}
where $q$ is the geometry order of approximation and we assume that integrals on $\Omega_h$ can be computed accurately. We further assume that there is a mapping $\Phi_h:\Omega_h \rightarrow \Omega$ that allows to map the approximated domain onto the exact domain. This mapping is assumed to be a piecewise smooth bijection, and fulfils $\Phi_h(\Gamma_h)=\Gamma$ and
\begin{equation}\label{eqn.mapping.bounds}
  \nrm{\Phi_h - \id}_{\xLinfty(\Omega_h)}\lesssim h^{q+1},\qquad
  \nrm{D\Phi_h - I}_{\xLinfty(\Omega_h)}\lesssim h^{q}.
\end{equation}

In this setting, we introduce adjusted versions of the previous bi- and linear forms and discrete norms of the DG discretisations. To this end, we effectively replace $\Omega$ by its approximation $\Omega_h$ yielding slightly modified discrete regions $\mathcal{T}_h$, $\mathcal{T}_h^\Gamma$, $\Omega_{\mathcal{T}}$\footnote{In regard to the discrete regions, we just assume $\mathcal{T}_h$, $\mathcal{T}_h^\Gamma$, $\Omega_{\mathcal{T}}$ to be defined accordingly from now on and---in the interest of readability---refrain from introducing new symbols.}, the forms $\Bhg(\cdot,\cdot)$, $\Ahg(\cdot,\cdot)$, $\Shg(\cdot,\cdot)$, $\Lhg(\cdot)$ and norms $\tnrmBg{\cdot}$, $\tnrmAg{\cdot}$, $\snrmSg{\cdot}$. Let us stress that all previous dependencies on $\Omega$, e.g., the selection of active elements $\Th$, are now replaced by dependencies on the discrete domain $\Omega_h$. As in Section~\ref{sec.error-exact-geometry}, we restrict ourselves to the homogeneous case $f=0$.
For the boundary data given by $\ghom\in \xHn{{1,\infty}}(\Gamma)$, we assume that there exists a sufficiently smooth extension into a domain $\Gamma^e\supset \Gamma\cup\Gamma_h$. We refer to the extension by $\ghom$, abusing the notation.
The plain-DG discretisation with geometry errors then reads as: Find $u_h \in \PP^k(\Th)$, such that
\begin{gather}\tag{geoDG\textsuperscript{hom}}\label{eqn.poisson-poblem.dg.wgeom}
  \Bhg(u_h, v_h) \coloneqq \Ahg(u_h,v_h) + \Shg(u_h,v_h) = \Lhg(0, \ghom; v_h) \quad \forall v_h \in \PP^k(\Th).
\end{gather}
Similarly, the Trefftz discretisation reads: Find $\uT \in \TT^k(\Th)$, such that
\begin{gather}\tag{geoTDG\textsuperscript{hom}}\label{eqn.poisson-poblem.trefftz.wgeom}
  \Bhg(\uT, v_h) \coloneqq \Ahg(\uT,v_h) + \Shg(\uT,v_h) = \Lhg(0, \ghom; v_h)\quad \forall v_h \in \TT^k(\Th).
\end{gather}

\subsubsection{A Priori Error Bounds}
For the error analysis, we introduce the auxiliary bilinear form with respect to the exact geometry. For $u, v\in \xHtwo(\Omega)\oplus\{v\cPhii\,\vert\,v\in V_h)$, we define
\begin{equation*}
  \A(u,v) \coloneqq (\nabla u, \nabla v)_{\Omega} - (\mean{\dnF u }, \jump{v})_{\Phi(\Fh) \cap \Omega} - (\dn u, v)_\Gamma,\quad\text{and}\quad \L(f;v) = (f,v)_{\Omega}.
\end{equation*}

We split the bilinear form $A_h$ into a part containing the symmetry and boundary control terms and another containing the remainder, i.e.\@{} the bulk, consistency and DG terms
\begin{align*}
  A_h^2(u,v) &= - (u, \dn v)_{\Gamma_h} + \beta (h^{-1} u,v)_{\Gamma_h} && A_h^1(u,v) = A_h(u,v) - A_h^2(u,v).
\end{align*}
We note that it is not necessary to split the linear form $L_h$ due to $f=0$.

We have the following Strang type lemma for the DG method:
\begin{lmm}\label{lem.strang.dg+tdg}
Let $u\in \xHtwo(\Omega)$ be the solution to \eqrefh{eqn.strong-poisson} with data $\ghom\in \xHn{{\frac{1}{2}+\varepsilon}}(\Omega\cup\Omega_h),~\varepsilon >0$. Furthermore, let either $V_h=\PP^k(\Th)$ and $u_h\in V_h$ be the solution to \eqref{eqn.poisson-poblem.dg.wgeom}, or let $V_h=\TT^k(\Th)$ and $u_h\in V_h$ be the solution to \eqref{eqn.poisson-poblem.trefftz.wgeom}. Then
  \begin{align}\label{eqn.strang-estimate.dg}
    \tnrmAg{u\cPhi - u_h} &\lesssim
        \inf_{v_h\in V_h}\left( \tnrmAg{u\cPhi - v_h} + \snrmSg{v_h}\right)\nonumber\\
        &\quad + \sup_{w_h \in V_h}\frac{\vert A_h^1(u \cPhi,w_h) - \A(u, w_h\cPhii)\vert}{\tnrmBg{w_h}}
+ h^{-1/2}\nrm{u \cPhi- \ghom}_{\Gamma_h}.
  \end{align}
\end{lmm}

\begin{proof}The proof follows the lines of~\cite{LR17, Leh17}. Let us denote $\usim = u \cPhi$. With the triangle inequality we have for arbitrary $v_h \in V_h$ that
\begin{equation*}
  \tnrmAg{\usim - u_h} \leq\tnrmAg{\usim - v_h} + \tnrmAg{v_h - u_h}.
\end{equation*}
With $w_h \coloneqq u_h - v_h$ and Lemma~\ref{lemma.poisson.dg}, we see
\begin{align*}
  \tnrmBg{w_h}^2 &\lesssim B_h(u_h - v_h, w_h) 
    =L_h(0,\ghom;w_h) - A_h(v_h, w_h) - S_h(v_h, w_h)\\
    &\lesssim \vert A_h (\usim -v_h, w_h) \vert + \vert L_h(0,\ghom;w_h) - A_h(\usim,w_h) \vert + \snrmSg{v_h}\snrmSg{w_h}\\
    &\lesssim \tnrmAg{\usim - v_h}\tnrmAg{w_h} + \vert A_h(\usim,w_h) - L_h(0,\ghom;w_h)\vert + \snrmSg{v_h}\snrmSg{w_h}
\end{align*}
Dividing by $\tnrmBg{w_h}$ and observing that both $\tnrmAg{\cdot}$ and $\snrmSg{\cdot}$ are dominated by $\tnrmBg{\cdot}$, we have
\begin{equation*}
  \tnrmAg{\usim - u_h} \lesssim \inf_{v_h\in V_h}\left( \tnrmAg{\usim - v_h} + \snrmS{v_h}\right) + \sup_{w_h \in V_h}\frac{\vert \Ahg(\usim,w_h) - \Lhg(0,\ghom;w_h)\vert}{\tnrmBg{w_h}}.
\end{equation*}
Due to integration by parts, we have $\A(u, w_h\cPhii) = \L(0; w_h\cPhii) = 0$ for all $w_h\in V_h$. Therefore,
\begin{align*}
  \vert\Ahg(\usim,w_h) - \Lhg(0,\ghom;w_h)\vert
    &= \vert \Ahg(\usim,w_h) - \Lhg(0,\ghom;w_h) - \A(u, w_h\cPhii)\vert\\
    &\leq \vert \Ahg^1(\usim,w_h) - \A(u, w_h\cPhii) \vert + \vert \Ahg^2(\usim,w_h) - \Lhg(0,\ghom;w_h)\vert.    
\end{align*}
For the final term, we observe
\begin{align*}
  \vert \Ahg^2(\usim,w_h) - \Lhg(0,\ghom;w_h)\vert
  &= \bigg\vert \int_{\Gamma_h}(-(\dn w_h )(\usim - \ghom) + w_h\frac{\beta}{h}(\usim - \ghom) \dif s \bigg\vert\\
  &\lesssim \tnrmBg{w_h}h^{-\frac{1}{2}}\nrm{\usim - \ghom}_{\Gamma_h}.
\end{align*}
\end{proof}

\begin{prpstn}[Unfitted DG error estimate.]\label{prop.unfitted-dg-err-er-with-geom}
  Let $u\in \xWn{{3,\infty}}(\Omega)$ if $k \leq 2$ or $u\in \xHn{{k+1}}(\Omega)$ if $k\geq 3$ be a solution to \eqrefh{eqn.strong-poisson}.
  Assume that $\ghom\in \xWn{{1,\infty}}(\Gamma)$ is extended sufficiently smooth into the discrete domain. 
  For the solution $u_h\in\PP^k$ to \eqref{eqn.poisson-poblem.dg.wgeom} there holds
  \begin{equation}
    \tnrmAg{u\cPhi - u_h} \lesssim (h^k + h^q )\left(S(u) + \nrm{\ghom}_{1,\infty,\Gamma}\right)
  \end{equation}
  with 
  \begin{equation}
    \nrm{\cdot}_{m,\infty,X}\coloneqq \max_{l\leq m}\nrm{D^l\cdot}_{\xLinfty(X)}, \quad \text{and} \quad    S(u) \coloneqq
    \begin{cases}
        \nrm{u}_{\xWn{{3,\infty}}(\Omega)} &\text{for } k \leq 2\\
        \nrm{u}_{\xHn{{k+1}}(\Omega)} &\text{for } k\geq3.
    \end{cases}
  \end{equation}
\end{prpstn}

\begin{proof}
We need to estimate the terms on the right-hand side of the Strang estimate in \eqref{eqn.strang-estimate.dg}. With a suitably continuous extension operator $\Ex\colon \xHn{{k+1}}(\Omega)\rightarrow \xHn{{k+1}}(\OT)$ (for $k \geq 3$) as above or $\Ex\colon \xWn{{3,\infty}}(\Omega)\rightarrow \xWn{{3,\infty}}(\OT)$ (for $k \geq 2$), we set $u^e\coloneqq \Ex u$. We then have
\begin{equation}\label{eqn.proof.err-est1}
  \tnrmAg{u\cPhi - v_h} \leq \tnrmAg{u\cPhi - u^e} + \tnrmAg{u^e - v_h}
\end{equation}
For the first term, we have
\begin{equation*}
   \tnrmAg{u\cPhi - u^e} \lesssim h^{q+\frac{1}{2}} S(u).
\end{equation*}
The proof of this bound follows along the lines of \cite[Lemma~12]{Leh17}. The additional interior-penalty facet contributions not treated in \cite[Lemma~12]{Leh17} have the same structure as the $\xLtwo(\Gamma_h)$-expressions treated in \cite[Lemma~12]{Leh17}. We note that the cases $k\leq2$ and $k\geq 3$ are distinguished as the proof requires $u\in \xWn{{2,\infty}}(\Omega)$
which is implied by $H^{\ell}(\Omega)$-regularity only for $\ell \geq 4$.
For the approximation part in \eqref{eqn.proof.err-est1}, we use the unfitted projection based on averaged Taylor/polynomials as above. With $v_h=\pi_h^k \Ex u =\pi_h^k u$ per construction, we have as above
\begin{equation}\label{eqn.proof.err-est2}
  \tnrmAg{u^e - v_h} = \tnrmAg{u^e - \pi_h^k u} \lesssim h^k \Vert u \Vert_{H^{k+1}(\Omega)} \lesssim h^k S(u).
\end{equation}
For the ghost-penalty semi-norm, we further have from the consistency error estimate \eqref{eqn.ghost-penalty-consitency} from above, that
\begin{equation*}
  \snrmS{\pi_h^k u} \lesssim h^k \Vert u \Vert_{H^{k+1}(\Omega)} \lesssim h^k S(u).
\end{equation*}

For consistency error contributions in the Strang estimate \eqref{eqn.strang-estimate.dg}, we have for $\ghom\in \xHn{{1,\infty}}(\Gamma)$ and assuming the data extension is bounded on $\nrm{\ghom}_{1,\infty,\Gamma}$, that
\begin{subequations}\label{eqn.strang-estimates}
\begin{align}
  \vert A_h^1(u \cPhi,w_h) - \A(u, w_h\cPhii)\vert &\lesssim h^q \nrm{u}_{\xHtwo(\Omega)}\tnrmAg{w_h}\label{eqn.strang-estimates.A}\\
h^{-1/2}\nrm{u \cPhi- \ghom}_{\Gamma_h} &\lesssim h^{q+\frac{1}{2}}\nrm{\ghom}_{1,\infty,\Gamma},\label{eqn.strang-estimates.bnd}
\end{align}
\end{subequations}
for $w_h\in V_h$. We provide the proofs in Appendix~\ref{appendix.proof.stang.geom.errs}.

Combining these estimates proves the claim.
\end{proof}

\begin{prpstn}
  Let $u\in \xWn{{3,\infty}}(\Omega)$ if $k \leq 2$ or $u\in \xHn{{k+1}}(\Omega)$ if $k\geq 3$ be a solution to \eqrefh{eqn.strong-poisson}. Assume that $\ghom\in \xHn{{1,\infty}}(\Gamma)$ is extended sufficiently smooth into the discrete domain. For the solution $u_\TT\in\TT^k$ to \eqref{eqn.poisson-poblem.trefftz.wgeom} there holds
  \begin{equation}
    \tnrmAg{u\cPhi - u_\TT} \lesssim (h^k + h^q )\left(S(u) + \nrm{\ghom}_{1,\infty,\Gamma}\right)
  \end{equation}
\end{prpstn}

\begin{proof}
The proof is an immediate consequence of the proof of Proposition~\ref{prop.unfitted-dg-err-er-with-geom}, by the observation that due to $\Delta u=0$ and the choice of our interpolation operator, we have in \eqref{eqn.proof.err-est2} that $v_h=\pi_h^k\uhom \in\TT^k$.
\end{proof}

\section{Implementational Aspects}
\label{sec.further-methods}

This section discusses some implementational aspects of the Trefftz and aggregated DG methods. In particular, we discuss how these approaches can be implemented in a general DG code without having to implement the basis functions for the Trefftz space or the DG space on general patches.

\subsection{Embedded Trefftz Method}\label{sec.embtrefftz}
In the previous section, we did not discuss the construction of $\TT^k$. One way is to set up the basis of harmonic polynomials; see, e.g., \cite{Her84}. A more flexible choice is to construct $\TT^k$ through an embedding in the corresponding DG space $\PP^k$. This has recently been introduced as the \emph{embedded Trefftz DG method} \cite{LS22}. One of the advantages of using the embedded Trefftz DG method in this setting is that it allows the easy construction of a particular solution so that we can also deal with the inhomogeneous problem. We sketch the approach in this section but note that there are no specific adjustments for the unfitted setting that needs to be taken compared to the embedding procedure introduced in \cite{LS22}.

Let $\{\phi_i\}_{i=1}^N$ be a basis of $\PP^k(\Th)$ and $\Giso:\RR^N\rightarrow\PP^k(\Th)$, $\Giso(\xbu) = \sum_{i=1}^N\xbu_i\phi_i$ be the Galerkin isomorphism mapping vectors in $\RR^N$ to finite element functions. Note that $\Giso(\eb_i) = \phi_i$ for the canonical unit vectors $\eb_i$. We then define the following matrices and vector for $i,j=1,\dots,N$
\begin{align*}
  (\Bb)_{ij} = \Bhg(\Giso(\eb_i), \Giso(\eb_j)) = \Bhg(\phi_i, \phi_j), &&
  (\lb)_i = \Lhg(\Giso(\eb_i)) = \Lhg(\phi_i), &&
  (\Wb)_{ij} = \langle \Delta \phi_i, \Delta \phi_j \rangle_{0,h}.
\end{align*}
Here, $\langle \cdot, \cdot \rangle_{0,h} \coloneqq \sum_{T\in\Th} \langle \cdot, \cdot \rangle_{T}$ is the element wise $\xLtwo$ -inner product on the active mesh. We observe that
\begin{equation*}
  \ker(\Delta) = \Giso(\ker(\Wb)).
\end{equation*}
As the Trefftz space is the kernel of $\Delta$ in $\PP^k(\Th)$, we are equivalently looking for a basis of $\ker(\Wb)$ to characterise the Trefftz space. As $\Wb$ is block diagonal, with the blocks corresponding to the elements of the active mesh, we construct $\Wb$ element wise.

For $T\in\Th$, let $\Wb_T\in\RR^{N_T\times N_T}$ be the block in $\Wb\in\RR^{N\times N}$ corresponding to the element $T$ with $N_T=\dim(\PP^k(T))$. The dimension of the kernel of $\Wb_T$ is $M_T= \dim(\ker(\Delta)) = N_T - L_T$ with $L_T = \dim(\range(\Delta))$ (e.g. in 2D $M_T = (k+2)(k+1)/ 2 - k (k - 1) / 2 = 2k + 1$). We can compute the kernel of $\Wb_T$ and collect the set of orthogonal basis vectors in the matrix $\Tb_T\in\RR^{N_k \times M_T}$. This can be done numerically, for example using a QR-decomposition or singular value
decomposition (SVD), see also Figure~\ref{fig.sketch.svd}. These blocks are then collected into a block matrix $\Tb\in\RR^{N \times M}$, with which we have the characterisation $\ker(\Wb) = \Tb\cdot\RR^M$.

\begin{figure}
  \begin{center}
        $\Wb_T =\!\! { \color{gray}
      \left(\begin{array}{@{}c@{~}c@{~}c@{~}c@{~}c@{~}c@{}} 
        \vertbar &   & \vertbar & \vertbar &  & \vertbar \\[0.6ex] 
        \ub_1 &  \!\ldots\! & \ub_{L_T} & \ub_{L_T+1} &\!\ldots\! & \ub_{N_T} 
        \\[0.6ex] \vertbar  &  & \vertbar&\vertbar&  & \vertbar 
      \end{array}\right)
      \!\!\cdot\!\!
      }
      \left(\begin{array}{@{~}c@{~}c@{~}c@{~}c@{~}c@{~}c@{~}} 
        \sigma_1 \\& \ddots \\&& \sigma_{L_T} \\&&& \orange{0} \\&&&& \orange{\ddots} \\&&&&& \orange{0}
      \end{array}\right)
      \!\!\cdot\!\!
      \left(\begin{array}{@{~}c@{~}c@{~}c@{~}} 
        \horzbar & \vb_1^T & \horzbar \\& \vdots & \\\horzbar & \vb_{L_T}^T & \horzbar \\\orange{\horzbar} & \orange{\vb_{L_T+1}^T} & \orange{\horzbar} \\& \orange{\vdots} & \\\orange{\horzbar} & \orange{\vb_{N_T}^T} & \orange{\horzbar}
      \end{array}\right)
      $ $\!\!\leadsto\!\!$ $ \Tb_T \!\!=\!\!$ \orange{$
      \left(\begin{array}{@{}c@{~}c@{~}c@{}} 
        \vertbar &   & \vertbar \\[0.6ex]  
        \vb_{L_T+1} &  \ldots & \vb_{N_T}
        \\[0.6ex] \vertbar  &  & \vertbar
      \end{array}\right)\!\!\!\!
    $}
  \end{center}
  \caption{Sketch of extraction of the kernel of the local matrix $\Wb_T$ based on a local SV decomposition.}
  \label{fig.sketch.svd}
\end{figure}

With the help of $\Tb$, we can now naturally define the Trefftz Galerkin isomorphism as $\GisoT:\RR^M\rightarrow\TT^k(\Th)$, $\GisoT(\xbu) = \Giso(\Tb\xbu)$. Having assembled $\Bb, \lb$ and $\Tb$ for the standard DG method, the system corresponding to \eqref{eqn.poisson-poblem.trefftz} reads as follows: Find $\ubT (=\GisoT^{-1}(\uT))$, such that
\begin{equation*}
  \Tb^T\Bb \Tb \ubT = \Tb^T\lb.
\end{equation*}
Note that the Trefftz system with matrices $\BbT = \Tb^T\Bb \Tb$ and $\lbT = \Tb^T\lb$ can also be assembled in an element-by-element fashion avoiding the setup of the (larger) matrices and vectors $\Bb$, $\lb$.

The embedded Trefftz approach allows the construction of a particular solution $u_f \in \PP^k(\Th)$ in a generic way. To this end, we compute element-wise $(\wb_T)_{i}=(f,\mathcal G(\eb_i))=(f,\Delta\phi_i)$ and define $(\ub_{f})_T=\Wb^\dagger_T \wb_T.$
Here $\Wb^\dagger_T$ denotes the pseudoinverse of the matrix $\Wb_T$, which can be obtained using the QR decomposition or SVD of the matrix $\Wb_T$, which may have already been computed when numerically computing the kernel of $\Wb_T$.
Then a particular solution is given by $u_f=\GisoT(\ub_f)$, and the Trefftz solution to the corresponding homogenised problem corresponds to the solution $\ubT$ to
\begin{equation*}
  \Tb^T\Bb \Tb \ubT = \Tb^T(\lb - \Bb \ub_f).
\end{equation*}

Let us summarise that the key idea of the embedded Trefftz approach is to exploit the characterisation of the Trefftz DG space as the kernel of differential operator on a standard DG space on a linear level. A similar approach can also be applied to implement the aggregated finite element method as we discuss next.

\subsection{Embedded Aggregated FEM}
In the aggregated finite element method, degrees of freedom on problematic elements, i.e. cut elements or ill-shaped cut elements, are marked as dependent degrees of freedom which are determined from degrees of freedom on uncut elements using an extrapolation. In~\cite{BVM18}, this interpolation is done geometrically based on a nodal representation of finite element degrees of freedom. Here, we want to discuss a different approach in the virtue of the embedded Trefftz DG approach. Since $\PP^k(\Thag)\subset \PP^k(\Th)$, we would like to characterise $\PP^k(\Thag)$ as the kernel of an operator in $\PP^k(\Th)$. In fact, we can identify the patch-wise jump operator, i.e.\@{} the facet-patch-local version of the ghost-penalty operator, as a corresponding suitable operator. 

Again, let $\{\phi_i\}_{i=1}^{N}$ be a basis of $\PP^k(\Th)$ and $\{\phi_i^T\}_{i=1}^{N_T}\subset \{\phi_i\}_{i=1}^{N}$ be the basis functions corresponding to an element $T\in\Th$. For a non-trivial patch $\Thw\subset\Th$, we define the ghost penalty operator facet-patch-wise as
\begin{equation*}
  s_{h,\omega}(v,w) = \sum_{F\in\FhTw} s_{h,F}(v,w)
\end{equation*}
For trivial patches $\Thw = \{T\}$, the ghost penalty operator is zero. For, $i,j=1,\dots,N$, we can define the matrix
\begin{equation*}
  (\Wb)_{ij} =   \Shg(\phi_j,\phi_i)
\end{equation*}
with local version $(\Wb_\omega)_{ij} =   s_{h,\omega}(\phi_j,\phi_i)$.
As with the weak Trefftz method, we have that
\begin{equation*}
  \ker(\Shg) = \Giso(\ker(\Wb)).
\end{equation*}
As before, the matrix $\Wb$ is block diagonal, with each non-trivial block corresponding to one patch and the zero blocks corresponding to all trivial patches. As before, we can therefore compute the blocks $\Wb_\omega$ independently, compute the kernel of dimension $N_\omega = \dim(\PP^k(\omega))$ numerically and collect the orthogonal basis vectors of $\ker(\Wb_\omega)$ in small matrices $\Tb_\omega$. Let $M=\dim(\PP^k(\Thag)) = \dim(\PP^k(\Th\setminus\ThC)) + N_T \vert \Ch \vert$. We then build a block diagonal matrix $\Tb\in\RR^{N\times M}$ from these blocks, together with an identity block corresponding to all trivial patches.

With the help of $\Tb$ we can now naturally define the aggregated Galerkin isomorphism as $\GagT:\RR^M\rightarrow\TT^k(\Th)$, $\GagT(\xbu) = \Giso(\Tb\xbu)$. 

With the notation as for the embedded Trefftz method, the embedded aggregate finite element method corresponds to finding $\ub_{ag}$ such that
\begin{equation*}
  \Tb^T\Bb \Tb \ub_{ag} = \Tb^T\lb.
\end{equation*}
Note that here the ghost penalty is part of $\Bb$, but is used to determine $\Tb$.

The generic nature of this formulation is an advantage over the extrapolation based approach in \cite{BVM18}. The extrapolation in \cite{BVM18} is based on a geometrical extrapolation exploiting a nodal representation of degrees of freedom. In contrast, our generic formulation allows for an implementation of the aggregation, which is independent of the choice of basis functions and can easily be generalised to curved elements.

\subsection{Embedded Aggregated Trefftz DG Methods}

Both previously discussed approaches to embed a special DG finite element space into the standard finite element space can be combined. We simply combine the previous components to obtain an embedding for the aggregated Trefftz DG method into the standard finite element space.
Defining 
\begin{equation*}
(\Wb)_{ij} =  \Shg(\phi_j,\phi_i) + \langle \Delta \phi_j, \Delta \phi_i \rangle_{0,h},
\end{equation*}
and determining the kernel of the corresponding operator in $\PP^k(\Th)$ yields the space of patch-wise harmonic polynomials, i.e.\@{} the aggregated Trefftz space. The procedure can be executed - as in the previous section - patch by patch and results in a linear system of the form
\begin{equation*}
\Tb^T\Bb \Tb \ubTag = \Tb^T(\lb - \Bb \ub_f)
\end{equation*}
with patch-wise block diagonal $\Tb$. Note that here the ghost penalty is again not part of $\Bb$, but is used to determine $\Tb$.

\section{Numerical Examples}
\label{sec.numerical-examples}

The previously described methods are implemented using \texttt{NGSTrefftz}
\cite{Sto22} and \texttt{ngsxfem} \cite{LHPvW21}, add-on packages to the finite
element library \texttt{NGSolve/Netgen} \cite{Sch14, Sch97}. The python scripts
implementing the methods discussed in this paper and the full numerical results
presented below are freely available in the zenodo
repository~\cite{HLSvW22_zenodo}. The stabilisation parameters for the
  interior penalty and the Nitsche method are fixed in all examples to $\beta =
  10 k^2$ and the ghost penalty scaling to $\gamma = 0.01$ if not mentioned otherwise.
  For the solution of linear systems arising in the numerical examples, we used sparse direct solvers.

\subsection{Example 1: Two Dimensions}
As a first example we consider a ring shaped geometry $\Omega = \{\xb\in\RR^2 \::\: 1 / 4 < \sqrt{\xb_1^2 +\xb_2^2} < 3 / 4 \}$. We take the exact solution to be the harmonic function $u=\exp(\xb_1) \sin(\xb_2)$ and set $g=u$. Note that as a result, the boundary data is \emph{exact on the approximated boundary}, consequently, no higher-order geometry approximation is necessary.

\begin{figure}
  \centering
  \includegraphics{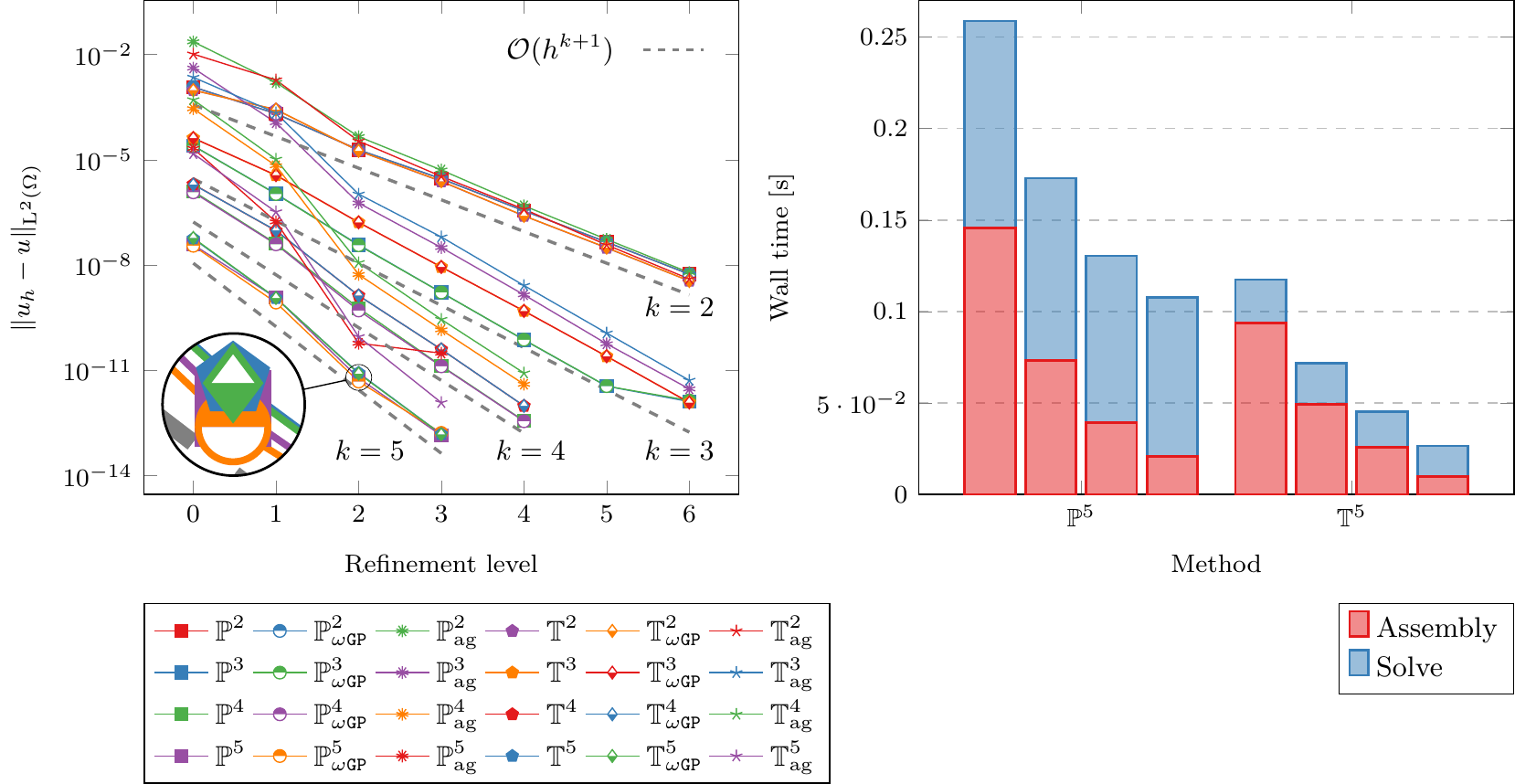}
  \caption{Example 1: Left: $\xLtwo$-error convergence for the unfitted Trefftz, embedded Trefftz and DG methods over a series of unstructured meshes with different polynomial orders. The exact boundary condition is applied on the approximated boundary and no higher-order geometry approximation is applied. Right: Compute times for specific steps of the methods with $h=0.1$, $k=5$ and using 1, 2, 4 and 12 parallel threads respectively.\protect\footnotemark}
  \label{fig.example1.l2err}
\end{figure}

From the methods analysed above, we consider the unfitted Trefftz method using both global ghost penalties, i.e., $\Fhgp=\Fhgpstar$, and patch-wise ghost penalties, i.e., $\Fhgp=\Fhgpmin$ and element aggregation.
The results using patch-wise ghost-penalty are indicated with a subscript ``${\omega\texttt{GP}}$'', the results using element aggregation using the subscript 
We also consider the embedded Trefftz method and the unfitted DG, both using global ghost-penalty stabilisation.
``ag'' and the embedded version of the Trefftz method is indicated with a subscript ``emb''.
The background domain is given by $\widetilde{\Omega}= (-1,1)^2$, the simplicial mesh is constructed by setting the mesh size $h=0.5^{i}$ for $i=1,\dots,7$, and we consider the orders $k=2,3,4,5$. We present the resulting $\xLtwo$-errors for the four methods in Figure~\ref{fig.example1.l2err}.

In Figure~\ref{fig.example1.l2err}, we see that all methods asymptotically converge with the optimal rate of $k+1$ in the $\xLtwo$-norm. We note that the DG methods appear to have a better error constant than the Trefftz methods. Furthermore, we see that the results between the patch-wise ghost penalties, global ghost penalty stabilisation choices are nearly indistinguishable. However, the element aggregation choice results in larger errors for both the Trefftz and full polynomial basis choice, and we observe pre-asymptotic behaviour for higher-order elements. This unsurprising, since element aggregation significantly reduces the number of elements in the case of very coarse meshes. Finally, the embedded Trefftz realisation of the Trefftz method lead to almost identical results. Consequently, they have been left out of the plot but are available in our archive \cite{HLSvW22_zenodo}. Regarding the performance of the different methods, the TDG method is indeed faster than the DG method, as expected, mainly due to the faster times to solve the resulting linear system. Furthermore, the exact compute times are given in Table~\ref{tab.timings.ex1} in Appendix~\ref{appendix:numerical_details}. The performance advantage of the embedded Trefftz method over the DG method is similar. For a detailed comparison between the Trefftz and embedded Trefftz methods, we refer to \cite{LS22}.

\footnotetext{Timings computed on an Intel Xenon E5-2687W v4 Processor @3.0GHz.}

\subsection{Example 2: Three Dimensions}
As a second example, we consider a three dimensional problem. The level set geometry is a flower-like shape taken from \cite{GM19}, and given by
\begin{equation*}
  \Omega = \{\xb\in\RR^3 : \sqrt{\xb_1^2 + \xb_2^2 + \xb_3^2} - 0.5 +  1 / 7 \cos(5 \atan_2(\xb_2, \xb_1)) \cos(\pi \xb_3) < 0\}.
\end{equation*}
We take the exact solution to be the harmonic function $u_{ex} =  \exp(\sqrt{2} \xb_1) \sin(\xb_2) \cos(\xb_3)$
and the boundary data is given by the exact solution. We apply the exact solution as the boundary data, and consequently the geometry error does not play a role here.

The background domain is $\Omega = (-1,1)^3$. The domain is meshed using structured tetrahedral elements with an initial mesh size of $h=0.5$, and we consider the polynomial order $k=2,3,4$. Furthermore, we only use the standard global ghost-penalty operator, as there was no visible difference between the two ghost penalty choices in the previous example.

\begin{figure}
  \centering
  \includegraphics{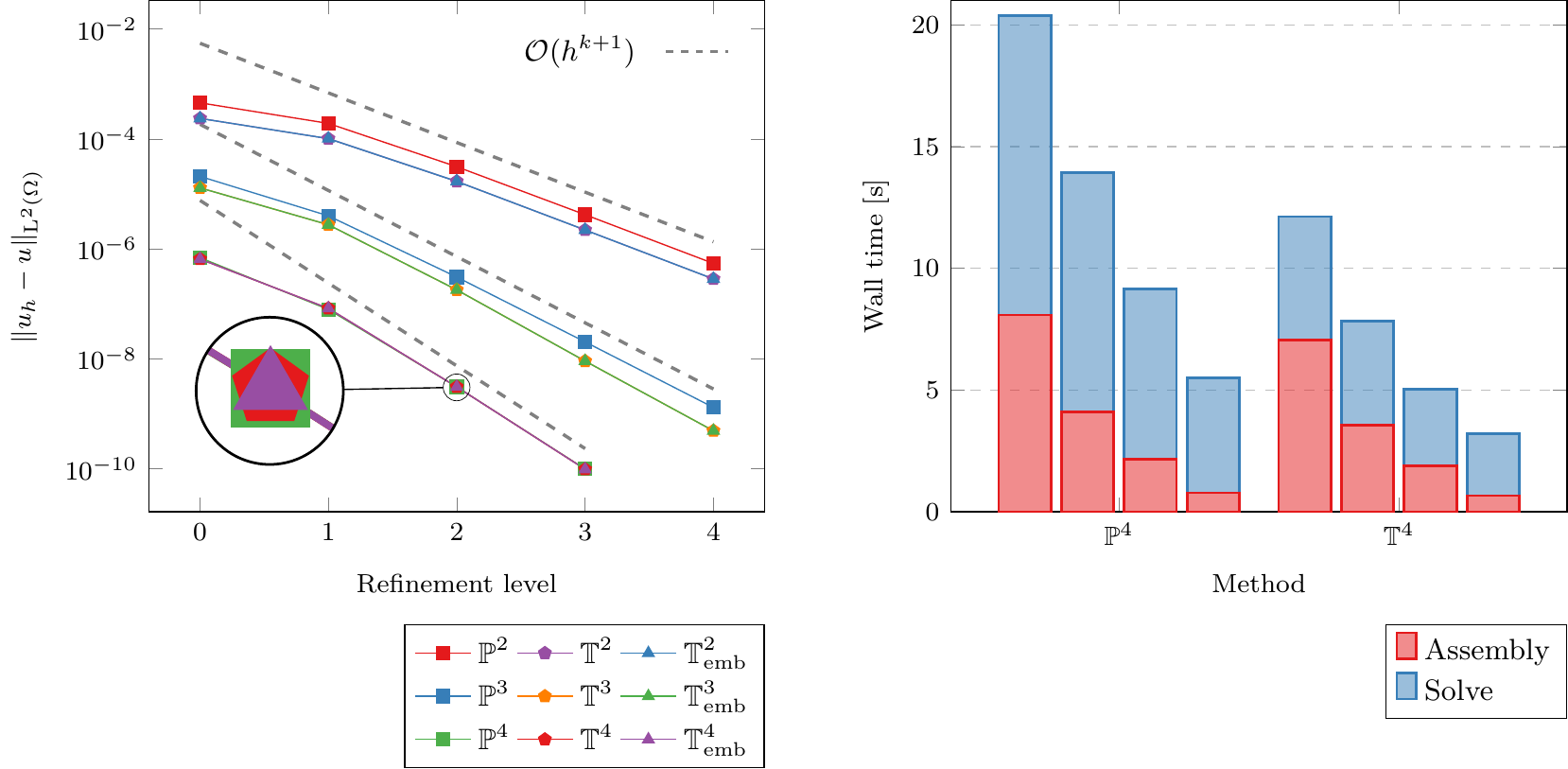}
  \caption{Example 2. Left: $\xLtwo$-error convergence for the unfitted Trefftz, embedded Trefftz and DG methods over a series of structured tetrahedral meshes with different polynomial orders. The exact boundary condition is applied on the approximated boundary and no higher-order geometry approximation is applied. Right: Compute times for specific steps of the methods with $h = 2^{-3}$, $k = 4$ and using 1, 2, 4 and 12 parallel threads respectively.\protect\footnotemark[4]}
  \label{fig.example3d_results}
\end{figure}

The results can be seen in Figure~\ref{fig.example3d_results}. We again observe optimal convergence in the $\xLtwo$-norm for all considered polynomial orders after some initial pre-asymptotic behaviour. In contrast to the previous examples, we see that the error constant favours the Trefftz methods rather than the DG method. Concerning performance, we again see that the TDG method is significantly faster than the DG method. The exact compute times are again presented in Appendix~\ref{appendix:numerical_details}.

\subsection{Example 3: Inhomogeneous Problem and Higher-order Geometry Approximation}
\label{sec.numerical-examples.subsec.ex3}
For our final example, we consider the same background domain and level set as in Example 1. With $r=\sqrt{\xb_1^2 + \xb_2^2}$, we choose the exact solution $u_{ex} = 20 (1 / 4 - r) (r - 3 / 4)$, which is zero on $\Gamma$. The right-hand side is taken as $f = -\Delta u_{ex}$.

\begin{figure}
  \begin{center}
    \includegraphics[width=0.5\textwidth]{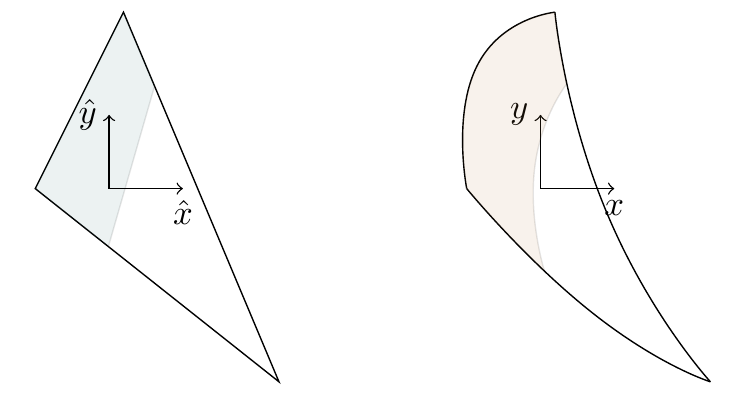}
  \end{center}
  \caption{Left: straight triangle and straight cut configuration, Right: curved element with higher order geometry approximation and coordinate system with respect to world coordinates}
  \label{fig.deormedmesh}
\end{figure}

As we can now apply the zero Dirichlet boundary condition on the discrete
interface, the geometry approximation error is present in this example. To
observe higher-order convergence, we, therefore, 
apply the strategy of isoparametric unfitted finite elements~\cite{Leh16,Leh17,LR17}. Here, a cheaply computable and small (magnitude $\lesssim h^2$) mesh deformation $\Theta_h$ is applied on the background mesh in a way such that a piecewise linear geometry approximation $\Omega^{\text{lin}}$ is mapped close to the exact geometry (in a higher order way), $\Omega_h = \Theta_h(\Omega^{\text{lin}}) \approx \Omega$. Thereby, the problem of numerical integration can be reformulated as the much simpler problem of numerical integration on a piecewise linear geometry. This enables robust higher-order convergence of the quadrature problem. 
At a first glance, changing the mesh may seem to contradict the unfitted finite element paradigm. However, we recall that the difficulty that is to be circumvented in unfitted finite element methods is the initial meshing problem (or the remeshing for moving domain problems) which is a non-local problem determining the mesh topology. In the isoparametric unfitted finite element approach, however, the necessary mesh alteration, the mesh deformation, is only local and does not change the mesh topology, i.e. the advantages of unfitted finite elements in terms of the meshing problem are not touched. We refer to \cite{Leh16} for further details.
In contrast to previous works with the isoparametric unfitted finite element methods, we exploit the flexibility of discontinuous Galerkin methods \footnote{In previous works the finite elements have been mapped according to the mesh deformation $\Theta_h$ to obtain ($H^1$-)conforming finite elements with, however, mapped polynomials $\tilde{p} = p \circ \Theta_h^{-1}$ where $p$ is a polynomial. As ($H^1$)-conformity is not needed with discontinuous Galerkin methods we don't need to apply a pull back to the reference geometry here.} and keep the coordinate system on the curved elements in world coordinates, see Figure~\ref{fig.deormedmesh}. 

We have seen in the previous two examples that the Trefftz and embedded Trefftz methods give indistinguishable results.
While the standard Trefftz method does not yield optimal approximation rates when applied to the inhomogeneous problem directly, the embedded Trefftz method immediately provides a particular solution to homogenise the problem with, as described in section~\ref{sec.embtrefftz}.
Therefore we only show results for the embedded Trefftz method here.
We again only use the global ghost-penalty choice.

\begin{figure}
  \centering
  \includegraphics{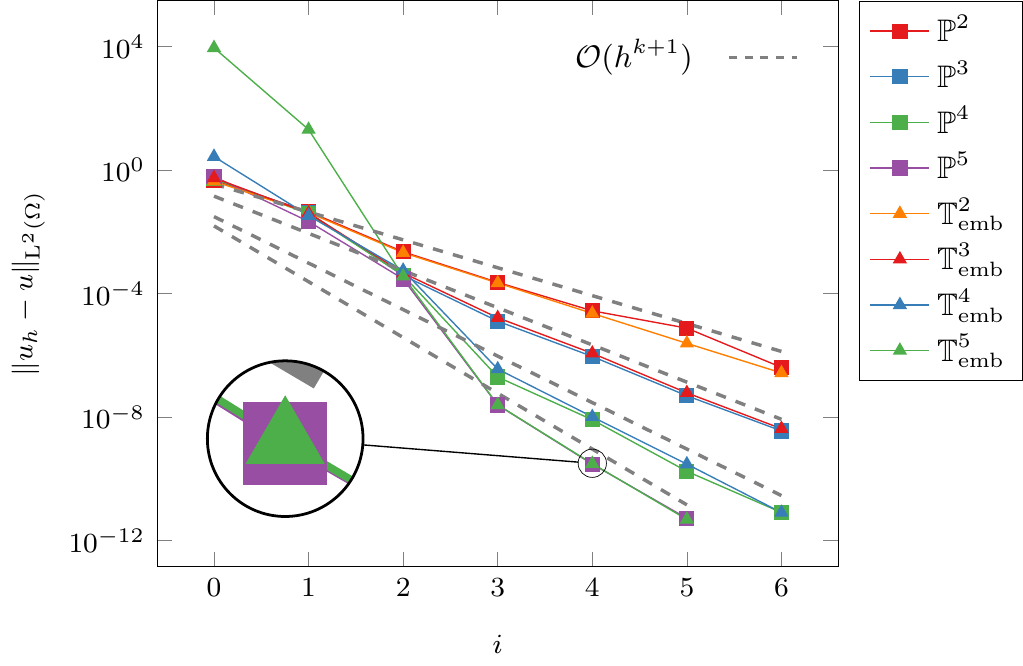}
  \caption{Example 3: $\xLtwo$-error convergence for the unfitted embedded Trefftz and DG methods over a series of meshes with different polynomial orders. The boundary is approximated to higher-order using a deformed mesh.}
  \label{fig.example3.l2err}
\end{figure}

The results can be seen in Figure~\ref{fig.example3.l2err}. Here we observe some pre-asymptotic behaviour on the three coarsest meshes. We attribute this to the fact that the width of the ring-shaped domain is of the same order of magnitude as the coarsest mesh size. As a result, the geometry and mesh deformation cannot be approximated properly on the coarser meshes. In particular, this becomes worse with increasing order of the deformation. However, we see optimal-order convergence for all orders and both methods once the mesh is sufficiently fine. We also see the slightly superior error constant for the DG method in the higher-order cases.

\subsection{Example 4: Species Dissolution From a Circle Subject to a Convection-Diffusion Equation}
\label{sec.numerical-examples.subsec.ex4}
This manuscript focuses on the Laplace equation as a model PDE to
introduce concepts for stable, reliable and arbitrarily high-order accurate
unfitted (Trefftz) DG methods and their analysis. In this last example, we
demonstrate -- without error analysis -- that the same concepts can be applied
to more general applications. Here, we consider a convection-dominated
convection-diffusion equation on a square background domain
$\widetilde{\Omega} = (-1, 1)^2$ with a circular obstacle $\Omega^\text{cir} =
B_R((0, 0))$ of radius $R = 1/4$, i.e., $\Omega = \widetilde{\Omega} \setminus
\Omega^\text{cir}$. A divergence-free flow field that is tangential on the
obstacle and essentially describes the advective transport from the left boundary
($x=-1$) to the right boundary ($x=1$) is given by $\bm{w} = (1 + R^2(y^2-x^2)/r^4, -2
R^2xy/r^4)^T$ and we set $c_{\mathbf{w}} = 2 \simeq \Vert \mathbf{w}
\Vert_{\infty}$ as a constant for the velocity magnitude. The convection-diffusion equation reads: 
Find $u\colon \Omega \to \mathbb{R}$ such that
\begin{subequations}\label{eqn.convection-diffusion.strong}
  \begin{align}
    - \alpha \Delta u + \bm{w} \cdot \nabla u &= 0 \quad \textnormal{in } \Omega, \\
    u &= 0 \quad \textnormal{on } \Gamma_l = \partial \widetilde{\Omega} \cap \{ x = -1\},\\
    \alpha \nabla u \cdot n &= 0 \quad \textnormal{on } \partial \widetilde{\Omega} \setminus \Gamma_l,\\
    u &= 1 \quad \textnormal{on } \partial \Omega^\text{cir}.
  \end{align}
\end{subequations}
We define $\Gamma \coloneqq \Gamma_l \cup \partial \Omega^{\text{cir}}$ as the part of
the boundary $\partial \Omega$ where Dirichlet boundary conditions are
prescribed ( so that the notation for the Nitsche formulation in \eqref{eq:Ah}
fits this setting).
We consider $\alpha = 10^{-3}$, i.e.\@{} a strongly convection-dominated
configuration with Péclet number $\operatorname{Pe} = \frac{c_{\mathbf{w}}
(2R)}{\alpha} = 1000$.
This problem is challenging for standard $H^1$-conforming Galerkin discretisations
due to a strong (parabolic) boundary layer that forms near the
obstacle, and proper convection stabilisation becomes necessary if the boundary
layer is not resolved, cf.\@{}, e.g., \cite{LR_SISC_2012}. In the context of unfitted DG methods considered here, we
can use a comparably simple upwind discretisation for the convection discretisation.

The discrete DG problem reads: Find $u_h \in \PP^k(\Th)$ such that
\begin{gather}\label{eqn.convection-diffusion-poblem.dg} \tag{CD-DG}
  \alpha \Ahg(u_h,v_h) + \Whg(u_h,v_h) + \Shg(u_h,v_h) = \Lhg(f,g;v_h) \quad \forall v_h \in \PP^k(\Th),
\end{gather}
where $\Whg(u,v)$ is the upwind DG bilinear form:
\begin{equation*}
  \begin{aligned}
    \Whg(u,v) \coloneqq &\ \sum_{T \in \Th} ( - u, \mathbf{w} \cdot \nabla v)_{T\cap\Omega} + ( \widehat{w_n u},  v)_{\partial (T \cap \Omega)} 
  \end{aligned}
\end{equation*}
with the upwind numerical flux $\widehat{w_n u} = \mathbf{w} \cdot
\mathbf{n}_{\partial T} \lim_{t \to 0^+} u(\mathbf{x} - t \mathbf{w})$ on
interior facets and outflow boundary facets and $\widehat{w_n u} = 0$ on all
inflow boundaries, i.e.\@{} boundaries with $\mathbf{w} \cdot
\mathbf{n}_{\Omega} < 0$.

For simplicity, we only consider the global ghost penalty stabilisation and
adjust the ghost penalty scaling to consider the contribution from the convection and choose
$\gamma = \gamma_0 (\alpha + h c_{\mathbf{w}})$ with $\gamma_0 = 0.001$. 

\begin{figure}
  \includegraphics[width=\textwidth,trim=13cm 54cm 13cm 0cm,clip=True]{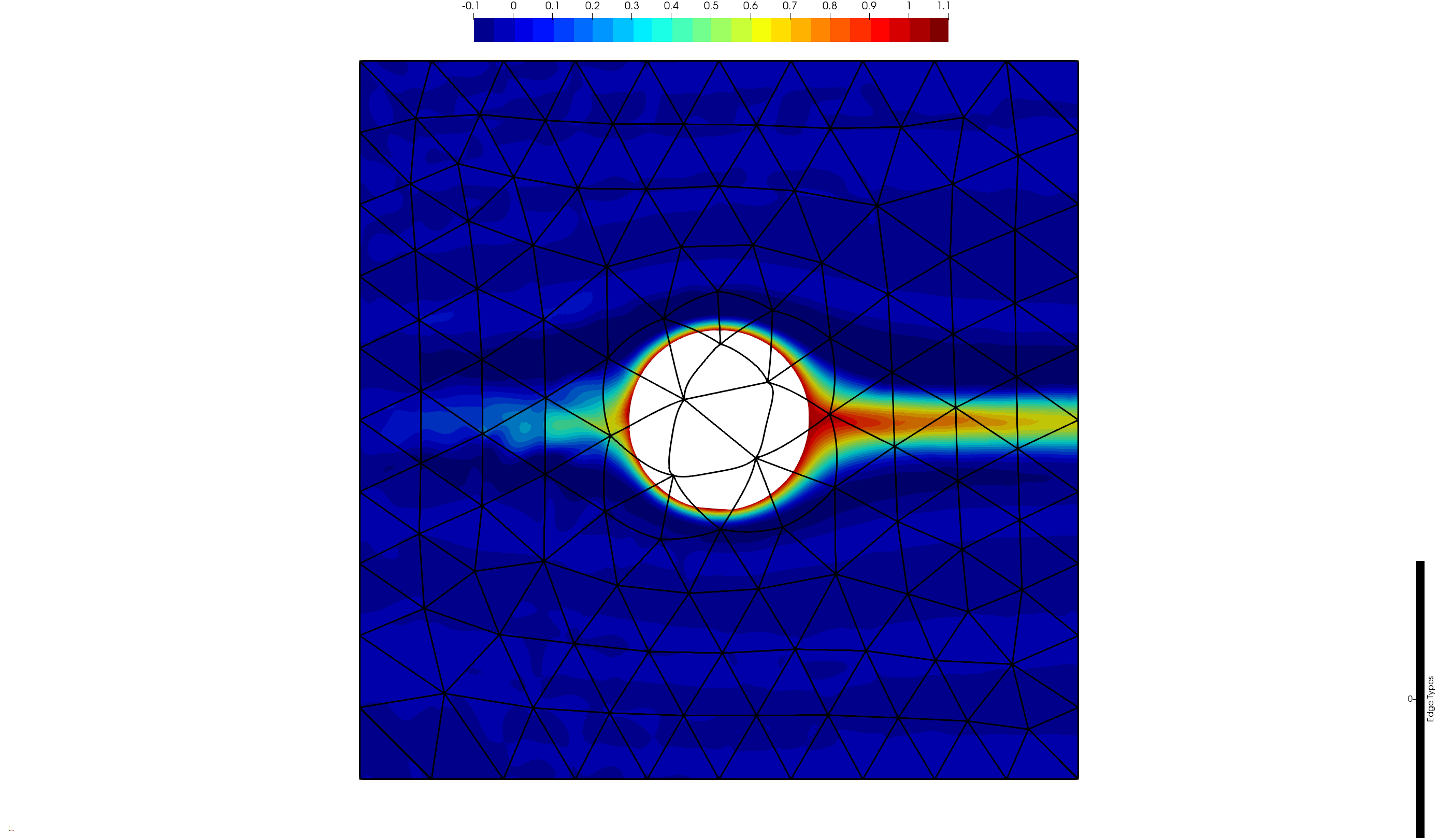} \\
  
  \centering
  \hspace*{-0.025\textwidth}
  \includegraphics[width=0.34\textwidth,trim=23.5cm 14cm 22.5cm 14cm,clip=True]{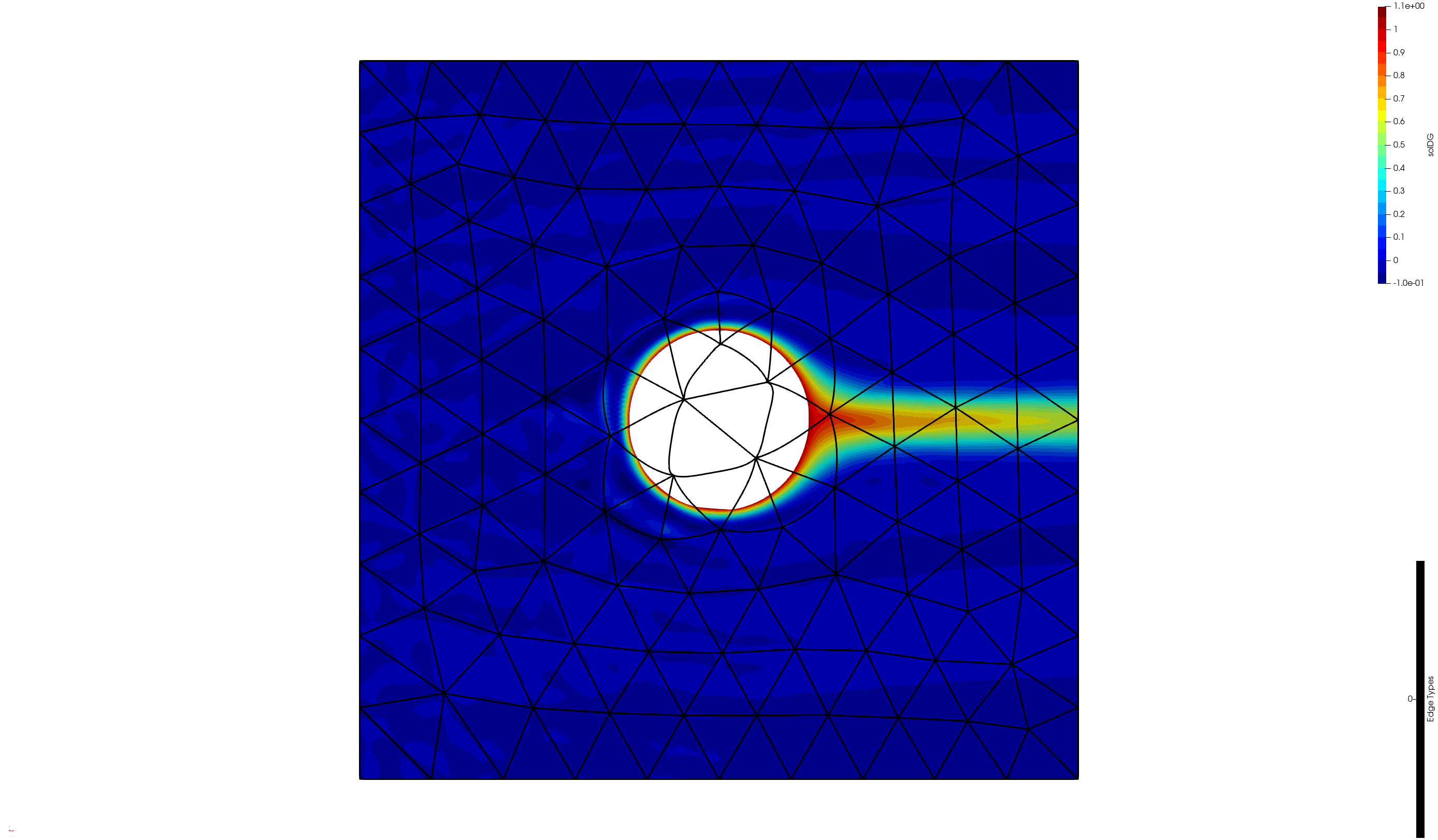}
  \hspace*{-0.025\textwidth}
  \includegraphics[width=0.34\textwidth,trim=23.5cm 14cm 22.5cm 14cm,clip=True]{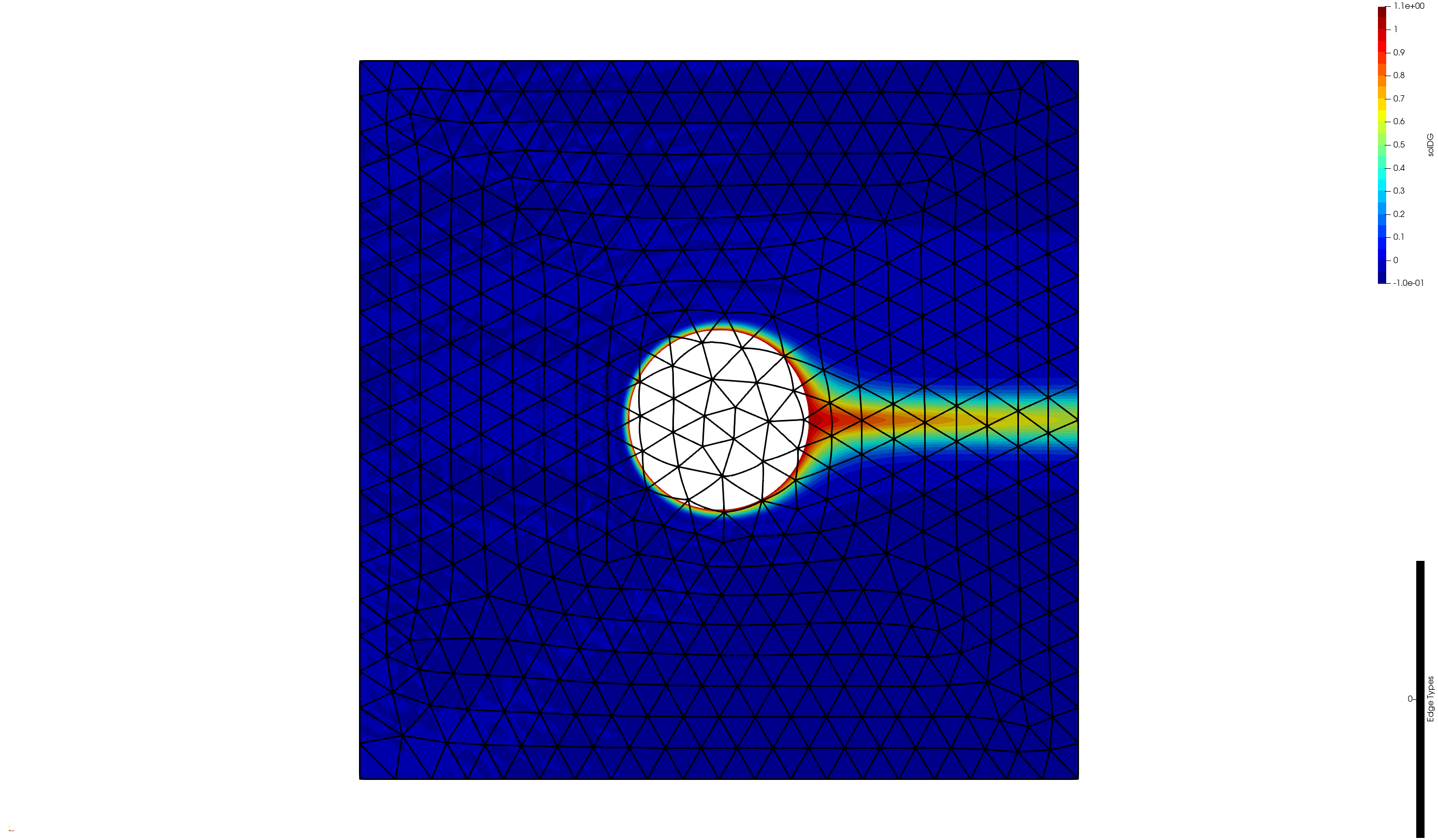}
  \hspace*{-0.025\textwidth}
  \includegraphics[width=0.34\textwidth,trim=23.5cm 14cm 22.5cm 14cm,clip=True]{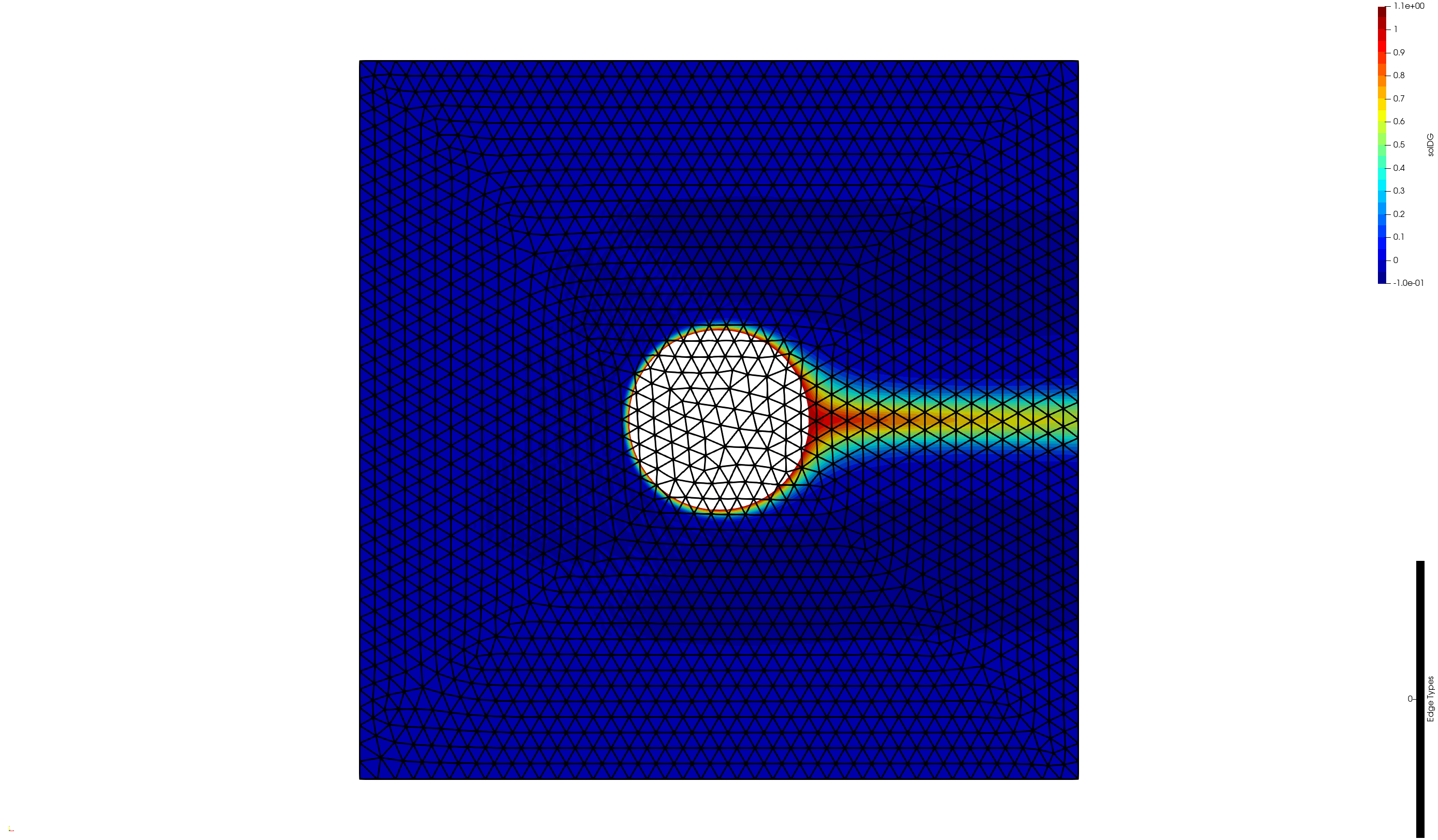}
  \hspace*{-0.025\textwidth}
  \\
  \hspace*{-0.025\textwidth}
  \includegraphics[width=0.34\textwidth,trim=23.5cm 14cm 22.5cm 14cm,clip=True]{plots/udg_cd_TDG0.png}
  \hspace*{-0.025\textwidth}
  \includegraphics[width=0.34\textwidth,trim=23.5cm 14cm 22.5cm 14cm,clip=True]{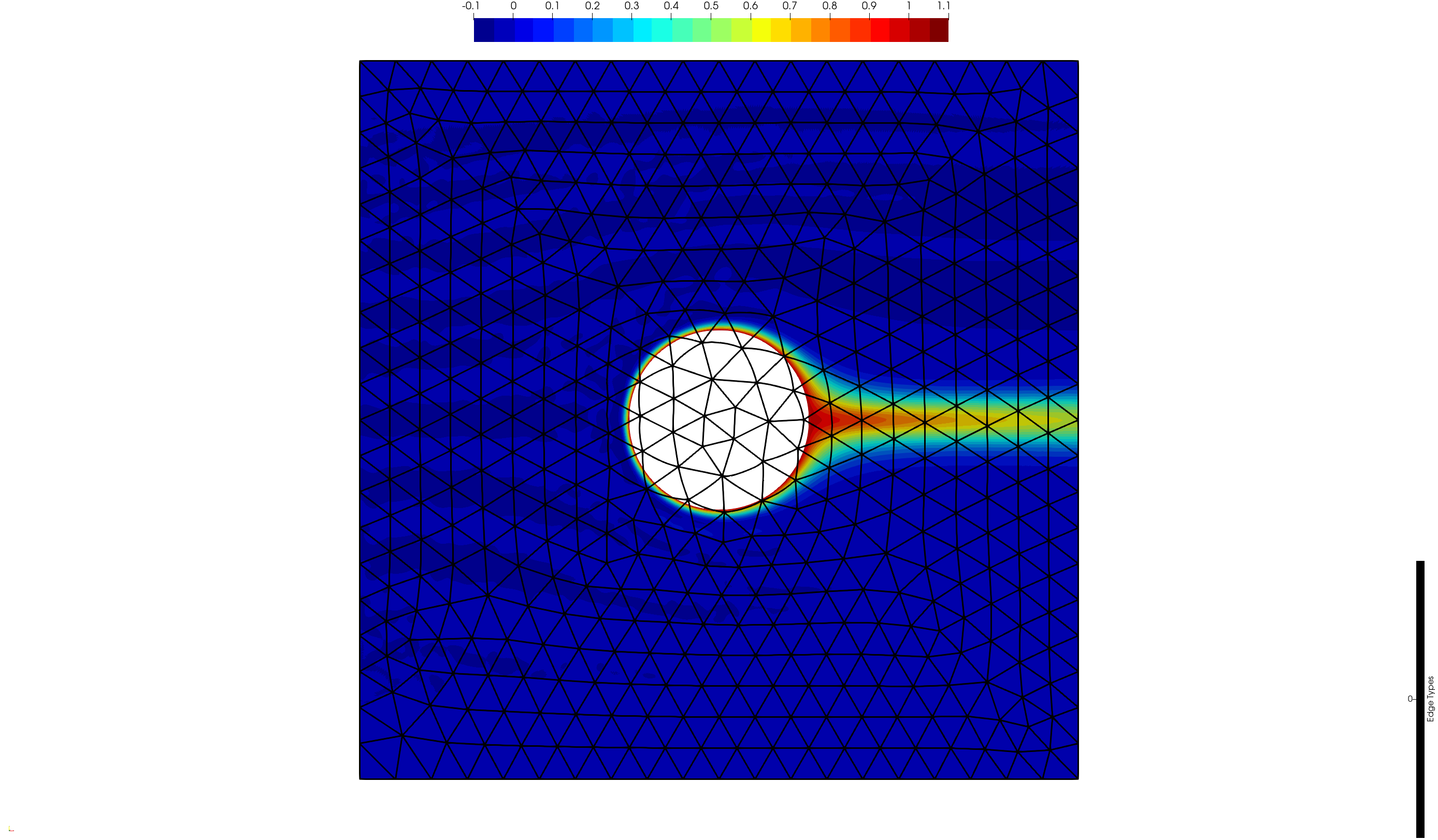}
  \hspace*{-0.025\textwidth}
  \includegraphics[width=0.34\textwidth,trim=23.5cm 14cm 22.5cm 14cm,clip=True]{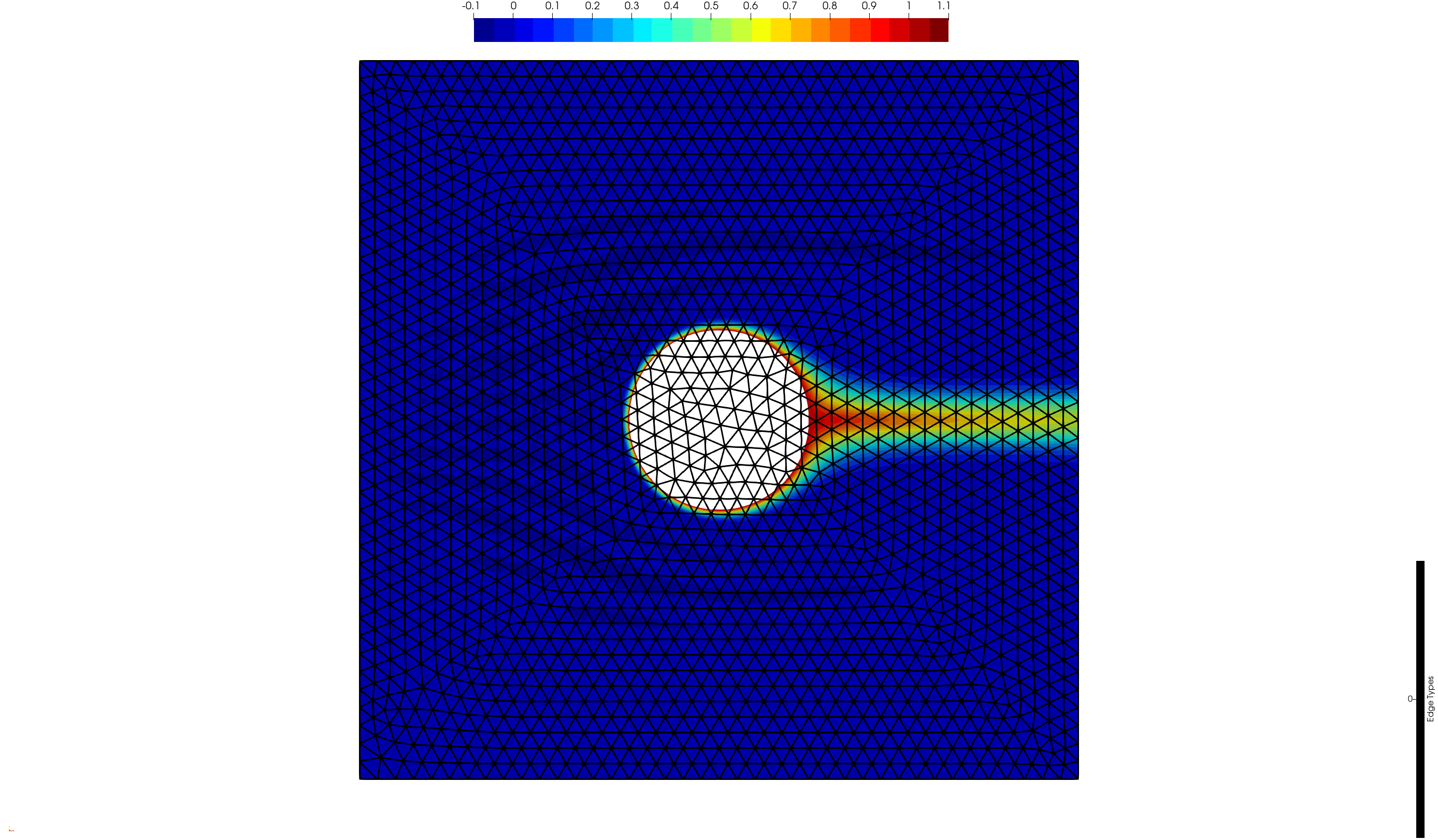}
  \hspace*{-0.025\textwidth}
\caption{Example 4: Top row: Unfitted DG solution to the convection-diffusion
    problem on curved meshes with $h=0.2, 0.1$ and $0.05$ for polynomial degree $k=4$, resulting in $3375,13530$ and $53565$ dofs respectively.
    Bottom row: Unfitted weak Trefftz DG solution with the same discretisation
    parameters and on the identical meshes resulting in $2025,8118$ and $32139$ dofs, respectively. The discrete colour scale emphasises
  small oscillations around zero.}
  \label{fig.example4}
\end{figure}

To define a proper Trefftz method, we can no longer use harmonic
polynomials but have to use a more generic construction of Trefftz basis
functions that are at least \emph{approximately} in the kernel of $\mathcal{L} = - \alpha
\Delta + \mathbf{w} \cdot \nabla$. To this end, we apply the idea of \emph{weak} Trefftz methods, implemented through an
embedding into the DG space as introduced in detail in \cite{LS22} and define
\begin{equation*}
  \TT^k(\Th) \coloneqq \{ v_h \in \PP^k(\Th) \mid \Pi_W (-\alpha \Delta v_h + \mathbf{w} \cdot\nabla
v_h )= 0\} ~\text{where}~ \Pi_W \text{ is the } L^2 \text { projection into } \PP^{k-2}(\Th).
\end{equation*}
Note that this is indeed a generalisation of the previously used Trefftz space
as we recover the space of harmonic polynomials for $\mathbf{w} = 0$ and $\alpha = 1$. Further, the thusly defined Trefftz DG space has the same dimension as the space of
harmonic polynomials, i.e.\@{} the same computational advantages over DG methods as
for the Laplace problem. The discrete Trefftz DG problem then reads: Find $u_h \in \TT^k(\Th)$ such that
\begin{gather}\label{eqn.convection-diffusion-poblem.tdg} \tag{CD-TDG}
  \alpha \Ahg(u_h,v_h) + \Whg(u_h,v_h) + \Shg(u_h,v_h) = \Lhg(f,g;v_h) \quad \forall v_h \in \TT^k(\Th).
\end{gather}
Figure~\ref{fig.example4} shows the results for the unfitted DG and the unfitted
Trefftz DG method of three decreasingly finer meshes (with $k=4$) where the meshes are curvilinear to obtain higher order geometry approximation (as in Section~\ref{sec.numerical-examples.subsec.ex3}). We observe that
strong oscillations are avoided by the upwind (and the ghost penalty)
stabilisation for both discretisations. Both schemes can capture the
characteristics of the problem even in these under-resolved situations with mesh
Péclet numbers  $\operatorname{Pe}_{h,k} = \frac{c_{\mathbf{w}}}{\alpha} \frac{h}{k} \in \{ 50, 25, 12.5\}$.
Corresponding to the characterisation of the weak Trefftz DG basis functions, we 
observe a slightly different pattern in the small oscillations around $0$ for
the Trefftz DG than for the DG method. While the high oscillations for the DG
method do not show an explicit direction, the small oscillations of the Trefftz
DG solution align with the flow field. However, the magnitude of the
oscillations for both methods are similar and decay under
refinement, as one would expect.

We have presented these results -- without going further
into the details -- to illustrate that the
methodology presented in this manuscript has considerable potential for a broader
class of unfitted PDE discretisations beyond the
Laplace equation.

\section*{Data availability statement}
The code used in this paper is available online in a Github repository: \url{https://github.com/hvonwah/unf-trefftz-poisson-code}, and is archived on zenodo \cite{HLSvW22_zenodo}.

\bibliographystyle{plain}
\bibliography{literature}

\begin{thebibliography}{10}

\bibitem{BNV22}
S.~Badia, E.~Neiva, and F.~Verdugo.
\newblock Linking ghost penalty and aggregated unfitted methods.
\newblock {\em Comput. Methods Appl. Mech. Engrg.}, 388:114232, 2022.

\bibitem{BVM18}
S.~Badia, F.~Verdugo, and A.~F. Mart{\'i}n.
\newblock The aggregated unfitted finite element method for elliptic problems.
\newblock {\em Comput. Methods Appl. Mech. Engrg.}, 336:533--553, 2018.

\bibitem{BE09}
P.~Bastian and C.~Engwer.
\newblock An unfitted finite element method using discontinuous {Galerkin}.
\newblock {\em Internat. J. Numer. Methods Engrg.}, 79(12):1557--1576,
  September 2009.

\bibitem{BMUP01}
T.~Belytschko, N.~Moes, S.~Usui, and C.~Parimi.
\newblock Arbitrary discontinuities in finite elements.
\newblock {\em Internat. J. Numer. Methods Engrg.}, 50:993--1013, 2001.

\bibitem{BS08}
S.~C. Brenner and L.~R. Scott.
\newblock {\em The mathematical theory of finite element methods}.
\newblock Springer, New York, 2008.

\bibitem{Bur10}
E.~Burman.
\newblock Ghost penalty.
\newblock {\em C.R. Math.}, 348(21-22):1217--1220, November 2010.

\bibitem{BCDE21}
E.~Burman, M.~Cicuttin, G.~Delay, and A.~Ern.
\newblock An unfitted hybrid high-order method with cell agglomeration for
  elliptic interface problems.
\newblock {\em SIAM J. Sci. Comput.}, 43(2):A859--A882, 2021.

\bibitem{B15}
E.~Burman, S.~Claus, P.~Hansbo, M.~G. Larson, and A.~Massing.
\newblock {CutFEM}: {D}iscretizing geometry and partial differential equations.
\newblock {\em Internat. J. Numer. Methods Engrg.}, 104:472--501, 2015.

\bibitem{BE18}
E.~Burman and A.~Ern.
\newblock An unfitted hybrid high-order method for elliptic interface problems.
\newblock {\em SIAM J. Numer. Anal.}, 56(3):1525--1546, 2018.

\bibitem{BE17}
E.~Burman and A.~Ern.
\newblock A cut cell hybrid high-order method for elliptic problems with curved
  boundaries.
\newblock In F.~Radu, K.~Kumar, I.~Berre, J.~Nordbotten, and I.~Pop, editors,
  {\em Numerical Mathematics and Advanced Applications ENUMATH 2017}, Lecture
  Notes in Computational Science and Engineering, pages 173--181, Cham, 2019.
  Springer.

\bibitem{BHLM17}
E.~Burman, P.~Hansbo, M.~G. Larson, and A.~Massing.
\newblock A cut discontinuous {Galerkin} method for the {Laplace}--{Beltrami}
  operator.
\newblock {\em IMA J. Numer. Anal.}, 37(1):138--169, 2017.

\bibitem{CDG21}
A.~Cangiani, Z.~Dong, and E.~Georgoulis.
\newblock hp-version discontinuous {Galerkin} methods on essentially
  arbitrarily-shaped elements.
\newblock {\em Math. Comp.}, 91(333):1--35, August 2021.

\bibitem{CDGH17}
A.~Cangiani, Z.~Dong, E.~H. Georgoulis, and P.~Houston.
\newblock {\em hp-Version Discontinuous Galerkin Methods on Polygonal and
  Polyhedral Meshes}.
\newblock Springer.

\bibitem{CD98}
O.~Cessenat and B.~Despres.
\newblock Application of an ultra weak variational formulation of elliptic
  {PDEs} to the two-dimensional {Helmholtz} problem.
\newblock {\em SIAM J. Numer. Anal.}, 35(1):255--299, 1998.

\bibitem{CGL09}
B.~Cockburn, J.~Gopalakrishnan, and R.~Lazarov.
\newblock Unified hybridization of discontinuous {Galerkin}, mixed, and
  continuous {Galerkin} methods for second order elliptic problems.
\newblock {\em SIAM J. Numer. Anal.}, 47(2):1319--1365, 2009.

\bibitem{DER14}
K.~Deckelnick, C.~M. Elliott, and T.~Ranner.
\newblock Unfitted finite element methods using bulk meshes for surface partial
  differential equations.
\newblock {\em SIAM J. Numer. Anal.}, 52(4):2137--2162, 2014.

\bibitem{ELL18}
D.~Elfverson, M.~G. Larson, and K.~Larsson.
\newblock {CutIGA} with basis function removal.
\newblock {\em Adv. Model. Simul. Eng. Sci.}, 5(1), March 2018.

\bibitem{ER12}
C.~M. Elliott and T.~Ranner.
\newblock Finite element analysis for a coupled bulk-surface partial
  differential equation.
\newblock {\em IMA J. Numer. Anal.}, 33(2):377--402, September 2012.

\bibitem{EH12}
C.~Engwer and F.~Heimann.
\newblock Dune-udg: {A} cut-cell framework for unfitted discontinuous
  {Galerkin} methods.
\newblock In {\em Advances in DUNE}, pages 89--100. Springer, 2012.

\bibitem{FOSS17}
T.~P. Fries, S.~Omerovi{\'{c}}, D.~Schöllhammer, and J.~Steidl.
\newblock Higher-order meshing of implicit geometries{\textemdash}part {I}:
  {I}ntegration and interpolation in cut elements.
\newblock {\em Comput. Methods Appl. Mech. Engrg.}, 313:759--784, January 2017.

\bibitem{GKF17}
C.~G\"{u}rkan, M.~Kronbichler, and S.~Fern\'{a}ndez-M\'{e}ndez.
\newblock {eX}tended hybridizable discontinuous {Galerkin} with heaviside
  enrichment for heat bimaterial problems.
\newblock {\em J. Sci. Comput.}, 72(2):542--567, August 2017.

\bibitem{GM19}
C.~G\"urkan and A.~Massing.
\newblock A stabilized cut discontinuous {Galerkin} framework for elliptic
  boundary value and interface problems.
\newblock {\em Comput. Methods Appl. Mech. Engrg.}, 348:466--499, May 2019.

\bibitem{HH02}
A.~Hansbo and P.~Hansbo.
\newblock An unfitted finite element method, based on {Nitsche}'s method, for
  elliptic interface problems.
\newblock {\em Comput. Methods Appl. Mech. Engrg.}, 191(47-48):5537--5552,
  November 2002.

\bibitem{HLSvW22_zenodo}
F.~Heimann, C.~Lehrenfeld, P.~Stocker, and H.~von Wahl.
\newblock Unfitted {Trefftz} discontinuous {Galerkin} methods for elliptic
  boundary value problems - {R}eproduction scripts, doi:
  \href{https://doi.org/10.5281/zenodo.8020304}{\texttt{10.5281/zenodo.8020304}},
  2022.

\bibitem{Her84}
I.~Herrera.
\newblock Trefftz method.
\newblock In {\em Topics in Boundary Element Research}, pages 225--253.
  Springer {US}, 1984.

\bibitem{HMP16}
R.~Hiptmair, A.~Moiola, and I.~Perugia.
\newblock Plane wave discontinuous {Galerkin} methods: exponential convergence
  of the {$hp$}-version.
\newblock {\em Found. Comput. Math.}, 16(3):637--675, 2016.

\bibitem{HMPS14}
R.~Hiptmair, A.~Moiola, I.~Perugia, and C.~Schwab.
\newblock Approximation by harmonic polynomials in star-shaped domains and
  exponential convergence of {T}refftz {$hp$-dGFEM}.
\newblock {\em ESAIM Math. Model. Numer. Anal.}, 48:727--752, May 2014.

\bibitem{HSK17}
S.~Hubrich, P.~Di Stolfo, L.~Kudela, S.~Kollmannsberger, E.~Rank, A.~Schröder,
  and A.~Düster.
\newblock Numerical integration of discontinuous functions: moment fitting and
  smart octree.
\newblock {\em Comput. Mech.}, 60(5):863--881, July 2017.

\bibitem{JL12}
A.~Johansson and M.~G. Larson.
\newblock A high order discontinuous {Galerkin} nitsche method for elliptic
  problems with fictitious boundary.
\newblock {\em Numer. Math.}, 123(4):607--628, September 2012.

\bibitem{K17}
F.~Kummer.
\newblock Extended discontinuous {Galerkin} methods for two-phase flows: {The}
  spatial discretization.
\newblock {\em Internat. J. Numer. Methods Engrg.}, 109(2):259--289, 2017.

\bibitem{LZ21}
M.~G. Larson and S.~Zahedi.
\newblock Conservative discontinuous cut finite element methods, May 2021.

\bibitem{Leh16}
C.~Lehrenfeld.
\newblock High order unfitted finite element methods on level set domains using
  isoparametric mappings.
\newblock {\em Comput. Methods Appl. Mech. Engrg.}, 300:716--733, March 2016.

\bibitem{Leh17}
C.~Lehrenfeld.
\newblock A higher order isoparametric fictitious domain method for level set
  domains.
\newblock In S.~Bordas, E.~Burman, M.~Larson, and M.~A. Olshanskii, editors,
  {\em Geometrically Unfitted Finite Element Methods and Applications -
  {P}roceedings of the UCL Workshop 2016}, volume 121 of {\em Lecture Notes in
  Computational Science and Engineering}, pages 65--92, Cham, 2017. Springer.

\bibitem{LHPvW21}
C.~Lehrenfeld, F.~Heimann, J.~Preu\ss, and H.~von Wahl.
\newblock \texttt{ngsxfem}: Add-on to {NGSolve} for geometrically unfitted
  finite element discretizations.
\newblock {\em J. Open Source Softw.}, 6(64):3237, August 2021.

\bibitem{LO19}
C.~Lehrenfeld and M.~A. Olshanskii.
\newblock An {Eulerian} finite element method for {PDE}s in time-dependent
  domains.
\newblock {\em ESAIM Math. Model. Numer. Anal.}, 53(2):585--614, March 2019.

\bibitem{LR_SISC_2012}
C.~Lehrenfeld and A.~Reusken.
\newblock {Nitsche-XFEM} with streamline diffusion stabilization for a
  two-phase mass transport problem.
\newblock {\em SIAM J. Sci. Comp.}, 34:2740--2759, 2012.

\bibitem{LR17}
C.~Lehrenfeld and A.~Reusken.
\newblock Analysis of a high-order unfitted finite element method for elliptic
  interface problems.
\newblock {\em IMA J. Numer. Anal.}, 38(3):1351--1387, August 2017.

\bibitem{LS22}
C.~Lehrenfeld and P.~Stocker.
\newblock Embedded {Trefftz} discontinuous {Galerkin} methods, January 2022.

\bibitem{LiShu2012}
F.~Li.
\newblock {On the negative-order norm accuracy of a local-structure-preserving
  {LDG} method}.
\newblock {\em {J. Sci. Comput.}}, 51(1):213--223, 2012.

\bibitem{LiShu2006}
F.~Li and C.-W. Shu.
\newblock {A local-structure-preserving local discontinuous {Galerkin} method
  for the {Laplace} equation}.
\newblock {\em {Methods Appl. Anal.}}, 13(2):215--234, 2006.

\bibitem{LX21}
S.~Lu and X.~Xu.
\newblock A geometrically consistent trace finite element method for the
  {Laplace}-{Beltrami} eigenvalue problem, 2021.

\bibitem{Mas12}
R.~Massjung.
\newblock An unfitted discontinuous {Galerkin} method applied to elliptic
  interface problems.
\newblock {\em SIAM J. Numer. Anal.}, 50(6):3134--3162, 2012.

\bibitem{MDB99}
N.~Moes, J.~Dolbow, and T.~Belytschko.
\newblock A finite element method for crack growth without remeshing.
\newblock {\em Internat. J. Numer. Methods Engrg.}, 46:131--150, 1999.

\bibitem{MKO13}
B.~Müller, F.~Kummer, and M.~Oberlack.
\newblock Highly accurate surface and volume integration on implicit domains by
  means of moment-fitting.
\newblock {\em Internat. J. Numer. Methods Engrg.}, 96(8):512--528, September
  2013.

\bibitem{OS16}
M.~A. Olshanskii and D.~Safin.
\newblock Numerical integration over implicitly defined domains for higher
  order unfitted finite element methods.
\newblock {\em Lobachevskii J. Math.}, 37(5):582--596, September 2016.

\bibitem{PDR07}
J.~Parvizian, A.~D\"uster, and E.~Rank.
\newblock Finite cell method.
\newblock {\em Comput. Mech.}, 41(1):121--133, 2007.

\bibitem{zbMATH01240671}
A.~Poullikkas, A.~Karageorghis, and G.~Georgiou.
\newblock The method of fundamental solutions for inhomogeneous elliptic
  problems.
\newblock {\em Comput. Mech.}, 22(1):100--107, 1998.

\bibitem{Pre18}
J.~Preu\ss.
\newblock Higher order unfitted isoparametric space-time {FEM} on moving
  domains.
\newblock Master's thesis, Georg-August-Universit\"at G\"ottingen, 2018.

\bibitem{Say15}
R.~I. Saye.
\newblock High-order quadrature methods for implicitly defined surfaces and
  volumes in hyperrectangles.
\newblock {\em SIAM J. Sci. Comput.}, 37(2):A993--A1019, January 2015.

\bibitem{Sch97}
J.~Sch\"oberl.
\newblock {NETGEN} an advancing front {2D}/{3D}-mesh generator based on
  abstract rules.
\newblock {\em Comput. Vis. Sci.}, 1(1):41--52, July 1997.

\bibitem{Sch14}
J.~Sch\"oberl.
\newblock C++11 implementation of finite elements in {NGSolve}.
\newblock Technical report, September 2014.

\bibitem{Ste70}
E.~M. Stein.
\newblock {\em Singular Integrals and Differentiability Properties of
  Functions}, volume~30 of {\em Princeton Mathematical Series}.
\newblock Princeton University Press, Princeton, NJ, 1970.

\bibitem{Sto22}
P.~Stocker.
\newblock \texttt{NGSTrefftz}: Add-on to {NGSolve} for {Trefftz} methods.
\newblock {\em J. Open Source Softw.}, 7(71):4135, March 2022.

\bibitem{trefftz1926}
E.~Trefftz.
\newblock Ein {G}egenst{\"u}ck zum {Ritzschen} {Verfahren}.
\newblock {\em Proc. 2nd Int. Cong. Appl. Mech., Zurich, 1926}, pages 131--137,
  1926.

\bibitem{zbMATH02139293}
A.~U\'sci{\l}owska-Gajda, J.~A. Ko{\l}odziej, M.~Cia{\l}kowski, and
  A.~Fr\k{a}ckowiak.
\newblock Comparison of two types of {Trefftz} method for the solution of
  inhomogeneous elliptic problems.
\newblock {\em {Comput. Assist. Mech. Eng. Sci.}}, 10(4):661--675, 2003.

\bibitem{yang2020trefftz}
J.~Yang, Mi. Potier-Ferry, K.~Akpama, H.~Hu, Y.~Koutsawa, H.~Tian, and D.~S.
  Z{\'e}z{\'e}.
\newblock Trefftz methods and {T}aylor series.
\newblock {\em Arch. Comput. Methods Eng.}, 27(3):673--690, 2020.

\end{thebibliography}

\appendix
\section{Appendix}
\label{sec.appendix}

\subsection{Proof of \texorpdfstring{Lemma~\ref{lem:traceineq2}}{Unfitted Trace Inequality}}\label{lem:traceineq2:proof}
\begin{proof}
  We apply a change of coordinates, the corresponding mapping is denoted by
  $\Phi_T$. We recall that we assumed $\Gamma \in \xCtwo$. Let $\bm{x}_T$ be a point in $T^\Gamma = T \cap \Gamma$ and $P_T$ be the tangential plane to $\bm{x}_T$. 
  First, we apply a translation with $- \bm{x}_T$ and scaling with $h_T$ to shrink $T$ to a domain of size $\mathcal{O}(1)$. Second, we apply a rotation, denoted by the orthogonal matrix $Q_T$, so that the tangential plane after transformation, i.e.\@{} $\Phi_T(P_T)$ aligns with the $x_1$-$x_2$ plane (or $x_1$ axis in 2D). We then have $\bm{\xi} = \Phi_T(\bm{x}) = h_T^{-1} Q_T \cdot (\bm{x} - \bm{x}_T)$ and denote $\hat{T} = \{ \Phi_T(x) \mid x \in T \}$. For sufficiently fine mesh sizes the unit outer normal of $\hat{T}^\Gamma = \Phi_T(T^\Gamma)$ is close to $\bm{e}_d$, the $d$th unit vector, 
especially there holds $(\bm{n},\bm{e}_d) \geq 1 - c_{\kappa} h \geq C \in (0,1)$ where $c_\kappa$ only depends on the maximum curvature of $T^\Gamma$ which is uniformly bounded.
  Furthermore we have that either $T \cap \Omega$ or $T \setminus \Omega$ is shape regular (and hence allows for the application of Lemma~\ref{lem:traceineq1}). We denote this part as $T^\ast$ and $\hat{T}^\ast = \Phi_T(T^\ast)$. Hence, we have for each $\hat v \in \xHn{1}(\hat{T})$ that there holds
  \begin{align*}
    \int_{\hat{T}^\ast} 2 \hat{v} \nabla \hat{v} \cdot \bm{e}_d dx = \int_{\hat{T}^\ast} \operatorname{div} (\hat{v}^2 \bm{e}_d) dx = \int_{\partial \hat{T}^\ast} \hat{v}^2 (\bm{n},\bm{e}_d) ds \geq (1 - c_{\kappa} h) \Vert \hat{v} \Vert_{\hat{T}^\Gamma}^2 - \Vert \hat{v} \Vert_{\partial \hat{T}}^2 \\ 
    \Longrightarrow\quad\Vert \hat{v} \Vert_{\hat{T}^\Gamma}^2 \lesssim \Vert \hat{v} \Vert_{\partial \hat{T}}^2 + \int_{\hat{T}^\ast} |\hat{v}| \Vert \nabla \hat{v} \Vert_2 dx
    \lesssim \Vert \hat{v} \Vert_{\hat{T}}^2 + \vert \hat{v} \vert_{\xHn{1}(\hat{T})}^2
  \end{align*}
  where in the last step we made use of Lemma~\ref{lem:traceineq1} and a Cauchy-Schwarz inequality. Applying a change of variables from $\hat T = \Phi_T(T)$ to $T$ yields
$$
\Vert v \Vert_{{T}^\Gamma}^2 \lesssim h_T^{-1} \Vert {v} \Vert_{{T}}^2 + h_T \vert {v} \vert_{\xHn{1}({T})}^2.
    $$
    \begin{figure}
      \begin{center}
        \includegraphics[scale=.55]{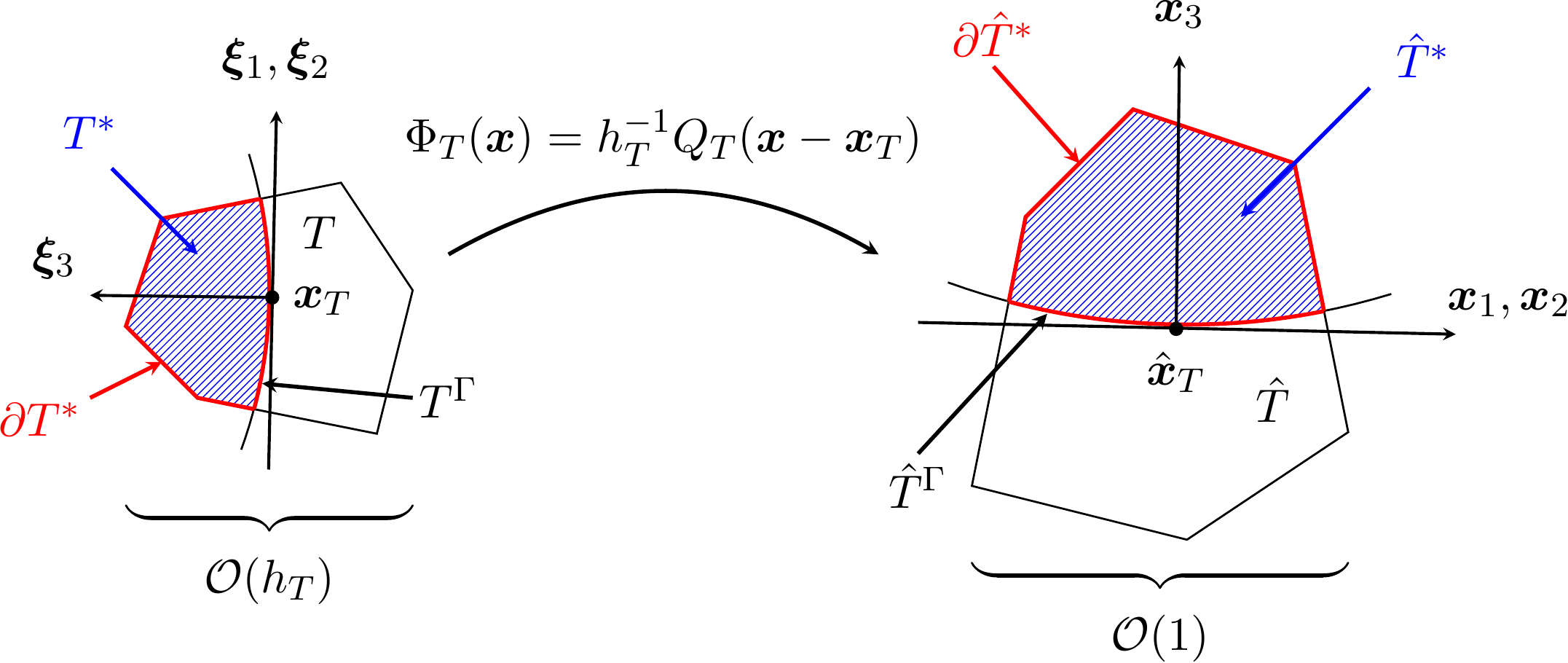}
      \end{center}
      \caption{Sketch of configuration in proof of Lemma~\ref{lem:traceineq2}.}
    \end{figure}

\end{proof}

\subsection{\texorpdfstring{Proof of Estimates \eqref{eqn.strang-estimates}}{Proof of Strang Geometry Error Terms}}
\label{appendix.proof.stang.geom.errs}

\begin{proof}[Proof of \eqref{eqn.strang-estimates.A}]
See also \cite[Appendix~A.4]{LR17}. For ease of notation let $\ut=u\cPhi$ and $\wht=w_h\cPhii$. We deal with the individual contribution in \eqref{eqn.strang-estimates.A} individually and begin with the volume terms. First we see using the chain-rule 
\begin{align*}
  \int_\Omega \nabla u \cdot \nabla\wht\dif\xbh = \int_{\Omega_h} \dDP\DPiT\nabla\ut\DPiT\nabla w_h\dif\xb\\
\end{align*}
This then gives with $J=\dDP$
\begin{multline}\label{eqn.strang-estimates.A-proof.volume}
\int_{\Omega_h} \nabla \ut \cdot \nabla w_h\dif\xbh - \int_\Omega \nabla u \cdot \nabla\wht\dif\xbh
 = \int_{\Omega_h} (\nabla \ut - J\DPiT\nabla\ut\DPiT)\cdot\nabla w_h\dif\xb\\
 \begin{aligned}[b]
  &=\int_{\Omega_h} \Big(\big((1 - J) + J(I - \DPiT)\big) \nabla \ut + J\DPiT\nabla\ut(I - \DPiT)\Big)\cdot\nabla w_h\dif\xb\\
  &\lesssim \bigg(\nrm{1 - J}_{\infty,\Omega_h}\nrm{\nabla\ut}_{\Omega_h} + \nrm{I - \DPiT}_{\infty,\Omega_h}\nrm{\nabla\ut}_{\Omega_h}
  + \nrm{\nabla\ut}_{\Omega_h}\nrm{I - \DPiT}_{\infty,\Omega_h} \bigg)\nrm{\nabla w_h}_{\Omega_h}\\
  &\lesssim h^q \nrm{\nabla u}_{\Omega}\nrm{\nabla w_h}_{\Omega_h}.
 \end{aligned}
\end{multline}
The final bound follows from \eqref{eqn.mapping.bounds} and the norm equivalence $\nrm{\ut}_{\Omega_h} \cong \nrm{u}_{\Omega}$ which can be proven along the same lines of argument.

For the symmetric interior penalty consistency term, let us consider single facet $F\in\Fh$. For simplicity, we identify $F$ with $F\cap\Omega$. Furthermore, we denote $\hat{F}=\Phi_h(F)$, $\nbF$ as the unit normal vector on $F$ and $\nbFh$ as the unit normal vector on $\hat{F}$. From \cite{LR17}, we then have
\begin{equation*}
  \int_{\hat{F}} \mean{\nbFh\cdot\nabla u }  \jump{ \wht } \dif\hat{s}
  = \int_{F} \DPiT \mean{\nabla \ut }\DPiT\nbF\jump{ w_h }\dif s.
\end{equation*}
It follows as above that
\begin{multline*}
 \int_{\hat{F}} \mean{\nbFh\cdot\nabla u } \jump{ \wht } \dif\hat{s} - \int_{F} \mean{\nbF\cdot\nabla \ut }\jump{ w_h } \dif s
 = \int_F \big(\dDP \DPiT \mean{\nabla \ut }\DPiT - \mean{\nabla \ut }\big)\nbF\jump{ w_h }\dif s\\
 \begin{aligned}[t]
  &= \int_F\big(C_1\mean{\nabla \ut } + C_2 \mean{\nabla \ut } C_3 \big) \nbF\jump{ w_h }\dif s\\
  &\lesssim (\nrm{C_1}_{\infty,\Omega_h} + \nrm{C_2}_{\infty,\Omega_h}\nrm{C_3}_{\infty,\Omega_h})\nrm{h^{1/2}\mean{\nabla\ut}\nbF}_{F}\nrm{h^{-1/2}\jump{w_h}}_{F}\\
  &\lesssim h^q\nrm{u}_{\xHtwo(\Omega)}\nrm{h^{-1/2}\jump{w_h}}_{F}
 \end{aligned}
\end{multline*}
with $C_1 = (1 - \dDP)I + \dDP(I - \DPiT)$, $C_2=\dDP\DPiT$ and $C_3=I - \DPiT$. The final estimate again follows from \eqref{eqn.mapping.bounds} and the continuos trace-estimate. Summing up over all facets then leads to the estimate
\begin{equation}\label{eqn.strang-estimates.A-proof.int}
 (\mean{\nbFh\cdot\nabla u }\jump{ \wht })_{\Phi_h(\Fh) \cap \Omega} - (\mean{\nbF\cdot\nabla\ut} \jump{ w_h })_{\Fh \cap \Omega} \lesssim h^q\nrm{u}_{\xHtwo(\Omega)} \tnrmAg{w_h}.
\end{equation}
The boundary consistency term estimate
\begin{equation}\label{eqn.strang-estimates.A-proof.bnd}
  (\nb_\Gamma\cdot\nabla u, \wht)_{\Gamma} - (\nb_{\Gamma_h}\cdot\nabla\ut, w_h)_{\Gamma_h}\lesssim h^q\nrm{u}_{\xHtwo(\Omega)} \tnrmAg{w_h}
\end{equation}
is proven analogously by replacing $F$ with $\Gamma_h$ and $\hat{F}$ with $\Gamma$.

For the symmetry and penalty terms of the symmetric interior penalty method, we observe that $\ut=u\cPhi$ is continuos on $\Omega_h$, and these terms vanish. Combing \eqref{eqn.strang-estimates.A-proof.volume} \eqref{eqn.strang-estimates.A-proof.int} and \eqref{eqn.strang-estimates.A-proof.bnd} proves the claim.
\end{proof}

\begin{proof}[Proof of \eqref{eqn.strang-estimates.bnd}]
See also \cite[Lemma~11]{Leh17}. We have that $u$ is equal to $g$ on $\Gamma=\Phi_h(\Gamma_h)$, and have assumed that $g$ is extended sufficiently smooth on $\Omega_h$ such that $\nrm{g^e}_{1,\infty,\Omega_h}\lesssim\nrm{g}_{1,\infty,\Gamma}$. Therefore,
\begin{equation}
  \nrm{u \cPhi- g}_{\Gamma_h}
  \lesssim \nrm{\Phi_h - \id}_{\infty,\Gamma_h}\nrm{g^e}_{1,\infty,\Omega_h} \lesssim h^{q+1}\nrm{g}_{1,\infty,\Gamma}.
\end{equation}
\end{proof}

\section{Additional Details on Numerical Results}
\label{appendix:numerical_details}

\begin{table}[!htb]
  \caption{Timing results for the DG and Trefftz DG method in Example 1 using polynomials of order 5 in two dimensions.}
  \label{tab.timings.ex1}
  \begin{tabular}{lrcll}
    \toprule
    Method & \# Dofs & \# Threads & Time Assemble [s] & Time Solve [s]\\
    \midrule
    DG  & 9051 & 1 & 0.14274429 & 0.1125111 \\
    TDG & 4741 & 1 & 0.09513908 & 0.0237321 \\
    \cmidrule(lr){3-5}
    DG  & 9051 & 2 & 0.07646305 & 0.0953993 \\
    TDG & 4741 & 2 & 0.04972788 & 0.0195056 \\
    \cmidrule(lr){3-5}
    DG  & 9051 & 4 & 0.04039365 & 0.0892400 \\
    TDG & 4741 & 4 & 0.02657581 & 0.0207315 \\
    \cmidrule(lr){3-5}
    DG  & 9051 & 12 & 0.02244950 & 0.1019655\\
    TDG & 4741 & 12 & 0.00980176 & 0.0154101\\
    \bottomrule  
  \end{tabular}
\end{table}

\begin{table}[!htb]
  \caption{Timing results for the DG and Trefftz DG method in Example 2 using polynomials of order 4 in three dimensions.}
  \label{tab.timings.ex2}
  \begin{tabular}{lrccc}
    \toprule
    Method & \# Dofs & \# Threads & Time Assemble [s] & Time Solve [s]\\
    \midrule
    DG  & 90370 & 1 & 8.13353458 & 12.1607583  \\
    TDG & 64550 & 1 & 7.09298664 &  5.0350659  \\
    \cmidrule(lr){3-5}
    DG  & 90370 & 2 & 4.12085618 & 10.6615596  \\
    TDG & 64550 & 2 & 3.68426016 &  4.3834508  \\
    \cmidrule(lr){3-5}
    DG  & 90370 & 4 & 2.15121641 & 6.6945196  \\
    TDG & 64550 & 4 & 1.86337521 & 2.8831929  \\
    \cmidrule(lr){3-5}
    DG  & 90370 & 12 & 0.74529837 & 4.55504362 \\
    TDG & 64550 & 12 & 0.66947819 & 2.81183520 \\
    \bottomrule  
  \end{tabular}
\end{table}

\end{document}